\pdfoutput=1
\documentclass[toc]{mathprint}
\usepackage{rutarmacros}
\usepackage{macros}
\title[Assouad spectra forms]{Attainable forms of Assouad spectra}

\author[Rutar]
  {Alex Rutar}
  {Mathematical Institute, University of St Andrews, St Andrews KY16 9SS, Scotland}
  {alex@rutar.org}

\begin{document}
\begin{abstract}
    Let $d\in\N$ and let $\varphi\colon(0,1)\to[0,d]$.
    We prove that there exists a set $F\subset\R^d$ such that $\dimAs\theta F=\varphi(\theta)$ for all $\theta\in(0,1)$ if and only if for every $0<\lambda<\theta<1$,
    \begin{equation*}
        0\leq (1-\lambda)\varphi(\lambda)-(1-\theta)\varphi(\theta)\leq (\theta-\lambda)\varphi\Bigl(\frac{\lambda}{\theta}\Bigr).
    \end{equation*}
    In particular, the following behaviours which have not previously been witnessed in any examples are possible: the Assouad spectrum can be non-monotonic on every open set, and can fail to be Hölder in a neighbourhood of 1.
\end{abstract}

\section{Introduction}
The Assouad dimension is a particular notion of dimension which captures the scaling properties of the ``thickest'' part of a set.
This in contrast to the more usual notions of box (or Hausdorff) dimension, which are in some sense a global measurement of scaling.
The Assouad dimension was originally introduced in \cite{ass1977} to study the embedding theory of metric spaces.
More recently, the Assouad dimension has received a significant amount of attention in the literature: see, for example, the books by Mackay \& Tyson on conformal geometry \cite{mt2010}, Robinson on embedding theory \cite{rob2010}, and Fraser on Assouad dimension in fractal geometry \cite{fra2021}.

If the box dimension and the Assouad dimension of a set agree, this implies that the set has a large amount of spatial regularity.
For instance, this is the case for any Ahlfors-regular subset of $\R^d$.
However, the box dimension and Assouad dimension can be distinct for many naturally-occurring sets, such as self-conformal sets with overlaps or self-affine sets.
In order to obtain a more fine-grained understanding of the Assouad dimension in this situation, the \defn{Assouad spectrum} was introduced by Fraser \& Yu in \cite{fy2018b}.
This is a notion of dimension parametrized by a variable $\theta\in(0,1)$, which approaches the box dimension as $\theta$ approaches 0 and the (quasi-)Assouad dimension as $\theta$ approaches 1.
We refer the reader to \cite{fra2021} for a general introduction to Assouad-type dimensions.

Besides being a useful bi-Lipschitz invariant and an important notion of fractal dimension in its own right, the Assouad spectrum provides more refined information about the Assouad dimension itself.
As a result, the Assouad spectrum has been explicitly studied for a wide range of examples (see, for example, \cite{bf2023,bff2022,fs2023,ft2021,fy2018a,fy2018b}).
This relationship has also been useful in applications outside of fractal geometry.
For instance, the Assouad spectrum plays an important role in the work by Roos \& Seeger \cite{zbl:07732556} on $L^p$ bounds for spherical maximal operators (and the analogous work in Heisenberg groups \cite{arxiv:2208.02774}).
The Assouad spectrum has also been used to obtain bounds for quasiconformal distortion in geometric mapping theory \cite{ct2023}.

In this paper, rather than consider explicit examples and applications of the Assouad spectrum, we focus on the general question of classification: what constraints on a function $\varphi\colon(0,1)\to[0,d]$ guarantee that there is a set $F\subset\R^d$ such that $\dimAs\theta F=\varphi(\theta)$ for all $\theta\in(0,1)$?

\subsection{Classifying Assouad spectra}
We fix $d\in\N$ and work in $\R^d$ with the Euclidean norm.
We write $B(x,\delta)$ to denote the open ball centred at $x$ with radius $\delta$.
If $F$ is a bounded subset of $\R^d$, for $\delta>0$, we let $N_\delta(F)$ denote the least number of balls of radius $\delta$ required to cover $F$.
Then, for $\theta\in(0,1)$, the \defn{Assouad spectrum} of $F\subset\R^d$ is given by
\begin{equation*}
    \dimAs\theta F=\inf\Bigl\{\alpha:(\exists C>0)(\forall 0 < \delta \leq 1)(\forall x\in F)\,\, N_{\delta^{1/\theta}}(F\cap B(x,\delta))\leq C\Bigl(\frac{\delta}{\delta^{1/\theta}}\Bigr)^\alpha\Bigr\}.
\end{equation*}
In general, $\lim_{\theta\to 0}\dimAs\theta F=\dimuB F$, and $\lim_{\theta\to 1}\dimAs\theta F=\dimqA F$ \cite{fhh+2019}, where $\dimqA F$ denotes the quasi-Assouad dimension of $F$ as introduced by Lü \& Xi \cite{lx2016}.
Like the Assouad dimension, the Assouad spectrum measures the worst-case local scaling of the set, but the Assouad spectrum specifies the relationship between the small and large scales.

The main result of this paper is to give a complete classification of possible forms of Assouad spectra.
\begin{itheorem}\label{it:assouad-char}
    Let $d\in\N$ and let $\varphi\colon(0,1)\to[0,d]$ be a function.
    Then there exists $F\subset\R^d$ such that $\dimAs\theta F=\varphi(\theta)$ for all $\theta\in(0,1)$ if and only if for every $0<\lambda<\theta<1$,
    \begin{equation}\label{e:assouad-upper-bd-intro}
        0\leq (1-\lambda)\varphi(\lambda)-(1-\theta)\varphi(\theta)\leq (\theta-\lambda)\varphi\Bigl(\frac{\lambda}{\theta}\Bigr).
    \end{equation}
\end{itheorem}
The proof of this result is given in \cref{s:class}.
For a geometric interpretation of the bound \cref{e:assouad-upper-bd-intro}, we refer the reader to \cref{ss:gint}.
The forward implication is well-known (see, for example, \cite[Theorem~3.3.1]{fra2021}); the reverse implication is proven in \cref{t:moran-exist}.

We can interpret the first inequality in \cref{e:assouad-upper-bd-intro} as a growth rate constraint, and the second inequality as an oscillation constraint.
In fact, the second inequality is always satisfied when $\varphi$ is increasing (the short argument is given in \cref{l:incr-bound}), which yields the following corollary.
\begin{icorollary}
    Let $d\in\N$ and let $\varphi\colon(0,1)\to[0,d]$ be an increasing function.
    Then there exists a set $F\subset\R^d$ with $\dimAs\theta F=\varphi(\theta)$ if and only if $\theta\mapsto(1-\theta)\varphi(\theta)$ is decreasing.
\end{icorollary}
We can also obtain results for the \defn{upper Assouad spectrum}, which is defined by bounding the lower scale from above, rather than specifying the relationship precisely:
\begin{align*}
    \dimuAs\theta F=\inf\Bigl\{\alpha:(\exists C>0)&(\forall 0 < \delta \leq 1)(\forall 0<\delta'\leq\delta^{1/\theta})(\forall x\in F)\\
                                                  & N_{\delta'}(B(x,\delta))\leq C\Bigl(\frac{\delta}{\delta'}\Bigr)^\alpha\Bigr\}.
\end{align*}
The upper Assouad spectrum is closely related to the Assouad spectrum: in fact $\dimuAs\theta F=\sup_{0<\theta'<\theta}\dimAs\theta F$ by \cite[Theorem~2.1]{fhh+2019}.
Combining this with \cref{c:upper-a} gives a full characterization of the upper Assouad spectrum (the details are given in \cref{ss:ua-char}).
\begin{icorollary}\label{ic:ua-char}
    Let $d\in\N$ and let ${\varphi}\colon(0,1)\to[0,d]$ be an arbitrary function.
    Then the following are equivalent:
    \begin{enumerate}[nl,a]
        \item\label{im:ua-as-set} There exists a set $F\subset\R^d$ such that $\dimuAs\theta F={\varphi}(\theta)$ for all $\theta\in(0,1)$.
        \item\label{im:ua-fun} ${\varphi}(\theta)$ is increasing and $\theta\mapsto(1-\theta){\varphi}(\theta)$ is decreasing.
        \item\label{im:ua-sup} ${\varphi}$ is the supremum of functions of the form $\theta\mapsto f(\theta)/(1-\theta)$ where
            \begin{equation*}
                f(\theta)=
                \begin{cases}
                    \kappa(1-c) &: 0<\theta\leq c\\
                    \kappa(1-\theta) &: c<\theta<1
                \end{cases}
            \end{equation*}
            for $c\in(0,1)$ and $\kappa\in[0,d]$.
    \end{enumerate}
\end{icorollary}
Beyond giving a full classification, \cref{it:assouad-char} also clarifies many of the properties of the Assouad spectrum: certain observations which might \emph{a priori} depend on explicit properties of the Assouad spectrum in fact only require the bound \cref{e:assouad-upper-bd-intro}.
For instance, the observation that if $\dimuB F=0$ then $\dimAs\theta F=0$ only requires the fact that $\lim_{\theta\to 0}\dimAs\theta F=\dimuB F$ along with the general bound (see \cref{p:a-props}).

We note that the 2-parameter family of functions in \cref{ic:ua-char} consists corresponds to the Assouad spectra of sets with upper box dimension $\kappa(1-c)$, quasi-Assouad dimension $\kappa$, and Assouad spectrum as large as possible.
In \cite{zbl:07732556}, such sets are called \defn{quasi-Assouad regular}.

Having completed the classification, in \cref{s:exceptional} we construct some exceptional sets.
Our first result concerns Hölder regularity.
\begin{itheorem}\label{it:non-holder}
    There is a compact set $F\subset\R$ such that $\theta\mapsto\dimAs\theta F$ is not Hölder in any neighbourhood of 1.
\end{itheorem}
In fact, there is no lower control on the rate at which $\dimAs\theta F$ approaches $\dimqA F$.
See \cref{ss:hol-failure} for the details.
This result is sharp: in \cref{p:lip-bounds}, we prove that $\dimAs\theta F$ is (uniformly) Lipschitz on $(0,1-\delta)$ for all $\delta>0$, with constants depending only on $\delta$ and the ambient dimension $d$.
This observation, along with \cref{it:non-holder}, provides a complete answer to \cite[Question~9.2]{fy2018b}.

Finally, we address the question of monotonicity.
In \cite[Question~17.7.1]{fra2021}, Fraser conjectures that the Assouad spectrum must be monotonic in some neighbourhood of $1$.
This was originally conjectured in \cite[Conjecture~2.4]{fhh+2019}.
We provide a strong negative answer to this question: we show that Assouad spectra that are non-monotonic on any open set are dense in the set of all possible upper Assouad spectra.
\begin{itheorem}\label{it:non-mono}
    For any $\varepsilon>0$ and function $\varphi$ satisfying one of the equivalent constraints in \cref{ic:ua-char}, there is a compact set $F\subset\R$ such that $\phi(\theta)=\dimAs\theta F$ is non-monotonic on any open subset of $(0,1)$ and $\norm{\phi-\varphi}_\infty<\varepsilon$.
\end{itheorem}
Since $\dimAs\theta F$ is Lipschitz on $(0,1-\delta)$ for every $\delta>0$, if $\varphi$ is non-constant then by Rademacher's theorem $\varphi$ must have strictly positive derivative on a set with positive Lebesgue measure.
This is sharp: using similar techniques as used in the proof of \cref{it:non-mono}, one can construct examples of sets with quasi-Assouad dimension $d$, box dimension arbitrarily close to 0, and Assouad spectrum that is strictly decreasing on a dense open subset of $(0,1)$ with Lebesgue measure arbitrarily close to 1.
We leave the details of such a construction to the interested reader.

\subsection{Rate constraints and the relationship with intermediate dimensions}
The intermediate dimensions are a different notion of dimension spectrum introduced in \cite{zbl:1448.28009} that interpolate between the Hausdorff and box dimensions.
In \cite{zbl:1509.28005}, the author and Banaji fully classify the possible forms of the intermediate dimensions.
For simplicity, in the discussion that follows we assume that the upper and lower intermediate dimensions coincide and denote the common value by $\dim_\theta F$.
We refer the reader to \cite{zbl:1509.28005} for precise statements of the results in full generality.

Recall that the (upper right) Dini derivative of a function $f\colon\R\to\R$ at $x$ is given by
\begin{equation}\label{e:Dini-intro}
    \diniu{+}f(x)=\limsup_{\varepsilon\searrow 0}\frac{f(x+\varepsilon)-f(x)}{\varepsilon}.
\end{equation}
We then recall the following result:
\begin{theorem}[\cite{zbl:1509.28005}]\label{it:intro-main-res}
    Let $g\colon(0,1)\to[0,d]$.
    Then there exists a non-empty bounded set $F\subset\R^d$ with $\dim_\theta F=g(\theta)$ if and only if
    \begin{equation}\label{e:intro-gen-bound}
        0\leq \diniu{+}g(\theta)\leq\frac{g(\theta)(d-g(\theta))}{d\cdot\theta}
    \end{equation}
    for all $\theta\in(0,1)$.
\end{theorem}
On the other hand, if $\varphi\colon(0,1)\to[0,d]$, \cref{ic:ua-char} gives that there exists $F\subset\R^d$ such that $\dimuAs\theta F=\varphi(\theta)$ if and only if
\begin{equation}\label{e:upper-growth-bd}
    0\leq \diniu+ \varphi(\theta)\leq \frac{\varphi(\theta)}{1-\theta}
\end{equation}
for all $\theta\in(0,1)$.
In particular, when $\dimB F<d$, $\diniu+g(\theta)\leq\frac{g(\theta)}{\theta}\cdot(d-\dimB F)/d$, so an arbitrary function which is the intermediate dimension of some set can be transformed to be the upper Assouad spectrum of a set through multiplication by a constant, reflection, and translation---and vice versa.

\subsection{Structure and outline of the paper}
In \cref{s:forms}, we study the family of functions $\mathcal{A}_d$ (see \cref{d:ad}) which satisfy the bound \cref{e:assouad-upper-bd-intro} for some fixed $d\in\N$.
This is the family which we will prove is the set of possible forms of Assouad spectra for subsets of $\R^d$.
First, in \cref{p:a-props}, we establish a number of basic properties of such functions.
The Assouad spectrum has been known to satisfy these properties, but here we only require the bound \cref{e:assouad-upper-bd-intro} and do not require any geometric properties of the Assouad spectrum itself.
Then in \cref{ss:rate}, we establish the growth rate bounds and the corresponding Lipschitz constraints.

Now, in \cref{s:class}, we prove \cref{it:assouad-char}.
The forward implication is standard, and follows by a straightforward covering argument: for completeness, we give the details in \cref{p:upper}.
To see the reverse implication, we will construct a homogeneous Moran set with prescribed Assouad spectrum using the techniques from \cite{zbl:1509.28005}.
This result is encapsulated in \cref{p:Moran-cnst}, where for a function satisfying certain derivative constraints, there exists a homogeneous Moran set such that the Assouad spectrum is given by a convenient formula.
It then remains to choose such a function carefully, which is done in \cref{t:moran-exist}.
In the remainder of the section, we construct families of monotonic and non-monotonic Assouad spectra in \cref{ss:non-mono-ex}, prove closure under suprema in \cref{ss:supr}, and complete the proof of \cref{ic:ua-char}.

To conclude, we use the classification result to construct examples of sets with exceptional Assouad spectra.
The result proving Hölder failure at 1 is described in \cref{ss:hol-failure}.
Then in \cref{ss:non-mono-open} we use the general family of non-monotonic spectra from \cref{ss:non-mono-ex} to construct a set with Assouad spectra which is not monotonic on any open subset of $(0,1)$.

\section{Forms of the family of functions \texorpdfstring{$\mathcal{A}_d$}{Ad}}\label{s:forms}
We first define the family of functions $\mathcal{A}_d$, which we will prove in \cref{s:class} are the possible forms of the maps $\theta\mapsto\dimAs\theta F$ for sets $F\subset\R^d$.
\begin{definition}\label{d:ad}
    Let $\mathcal{A}_d$ denote the set of functions $\varphi\colon(0,1)\to[0,d]$ where for any $0<\lambda<\theta<1$,
    \begin{equation}\label{e:assouad-upper-bd}
        0\leq (1-\lambda)\varphi(\lambda)-(1-\theta)\varphi(\theta)\leq(\theta-\lambda)\varphi\Bigl(\frac{\lambda}{\theta}\Bigr).
    \end{equation}
\end{definition}
In this section, we study properties of the family $\mathcal{A}_d$ directly: we emphasize that we do not require any geometric facts about the Assouad spectrum itself.

In \cref{p:a-props}, we will prove that functions in $\mathcal{A}_d$ are uniformly continuous.
Thus, we will embed $\mathcal{A}_d$ in $C([0,1])$ by defining $\varphi(0)=\lim_{\theta\to 0}\varphi(\theta)$ and $\varphi(1)=\lim_{\theta\to 1}\varphi(\theta)$.
We will use this notation once we prove uniform continuity.

\subsection{Rescaling and a geometric interpretation of the bound}\label{ss:gint}
Given $\varphi\in\mathcal{A}_d$, define $\beta(\theta)=(1-\theta)\varphi(\theta)$.
In \cref{e:assouad-upper-bd}, the first inequality implies that $\beta(\theta)$ is decreasing, and the second states that for all $0<\lambda<\theta<1$,
\begin{equation}\label{e:h-bd}
    \frac{\beta(\lambda)-\beta(\theta)}{\theta-\lambda}\leq\frac{\beta\bigl(\frac{\lambda}{\theta}\bigr)}{1-\frac{\lambda}{\theta}}.
\end{equation}
The left hand side is the negative of the slope of the line passing through $(\lambda,\beta(\lambda))$ and $(\theta,\beta(\theta))$, and the right hand side is the negative of the slope of the line passing through $(\lambda/\theta,\beta(\lambda/\theta))$ and $(1,0)$.
The secants in this constraint for a function $\beta$ are depicted in \cref{f:decr-bound}.
\begin{figure}[t]
    \begin{center}
        \begin{tikzpicture}[scale=8,font=\tiny,>=stealth]
            \draw[thick,->] (-0.05,0) -- (1.05,0);
            \draw[thick,->] (0,-0.05) -- (0,0.51825397);

            \matrix [draw,below left,scale=0.8,column sep=0.3em, inner sep=0.3em, align=left] at (current bounding box.north east) {
                \draw[thick] (0,-0.1) -- (0.71,-0.1); & \node[inner sep=0pt] {$\beta$}; \\
            };

            \begin{scope}[thick, smooth, variable=\x]
                \draw[domain=0:1/3] plot ({\x}, {(1-(\x)^2)-7/18-1/7});
                \draw[domain=1/3:2/3] plot ({\x}, {(\x-2/3)^2+7/18-1/7});
                \draw[domain=2/3:1] plot ({\x}, {-2.21428571*(2/3-\x)^2+7/18-1/7});
            \end{scope}
            \coordinate (lam) at (1/3,0.35714286);
            \coordinate (the) at (1/2,0.27380953);
            \coordinate (lat) at (2/3,0.24603175);
            \coordinate (one) at (1,0);

            \begin{scope}[thick,dashed]
                \draw (lam) -- (the);
                \draw (lat) -- (one);
            \end{scope}
            \draw (0.57735027,-0.01) node[below] {$\lambda^{1/2}$} -- (0.57735027,0.01);

            \begin{scope}[dotted]
                \draw (1/3,0) node[below]{$\lambda$}-- (lam);
                \draw (1/2,0) node[below]{$\theta$}-- (the);
                \draw (2/3,0) node[below]{$\lambda/\theta$}-- (lat);
                \node[below] at (one) {$1$};
            \end{scope}
        \end{tikzpicture}
    \end{center}
    \caption{A plot of $\beta(\theta)=(1-\theta)\varphi(\theta)$ where $\varphi\in\mathcal{A}_d$, and the lines with slopes corresponding to \cref{e:h-bd}.}
    \label{f:decr-bound}
\end{figure}

\subsection{Basic properties}
In this section, we collect various properties of the family $\mathcal{A}_d$.
First, we observe the following useful lemma which was essentially proven in \cite[Remark~3.8]{fy2018b}.
Here, we obtain it as a direct consequence of \cref{e:assouad-upper-bd}.
Heuristically, this lemma states that the function $\varphi(\theta)$ is ``almost increasing'', up to some possible local oscillations.
\begin{lemma}\label{l:nth-root-bd}
    Let $\varphi\in\mathcal{A}_d$.
    Given $0<\theta_1<\theta_2<\cdots<\theta_n<1$,
    \begin{equation*}
        \varphi(\theta_1)\leq \max\left\{\varphi\Bigl(\frac{\theta_1}{\theta_2}\Bigr),\varphi\Bigl(\frac{\theta_2}{\theta_3}\Bigr),\ldots,\varphi\Bigl(\frac{\theta_{n-1}}{\theta_n}\Bigr),\varphi(\theta_n)\right\}.
    \end{equation*}
    In particular, for any $n\in\N$ and $\theta\in(0,1)$, $\varphi(\theta)\leq \varphi(\theta^{1/n})$.
\end{lemma}
\begin{proof}
    Let $0<\theta_1<\theta_2<\cdots<\theta_n<1$.
    Applying \cref{e:assouad-upper-bd} to each pair $\theta_i,\theta_{i+1}$,
    \begin{equation*}
        (1-\theta_1)\varphi(\theta_1)\leq(1-\theta_n)\varphi(\theta_n)+\sum_{k=2}^n(\theta_k-\theta_{k-1})\varphi\Bigl(\frac{\theta_{k-1}}{\theta_{k}}\Bigr)
    \end{equation*}
    from which the result follows.
    Taking $\theta_i=\theta^{\frac{n-i+1}{n}}$ for each $i=1,\ldots,n$, observe that $\theta_{k-1}/\theta_k=\theta^{1/n}$ and $\theta_n=\theta^{1/n}$ so that $\varphi(\theta)\leq \varphi(\theta^{1/n})$.
\end{proof}
We now have the following essential properties of $\mathcal{A}_d$.
All of these properties have been previously observed for the Assouad spectrum, but the main point here is that these properties only depend on the family $\mathcal{A}_d$ and not on other properties of the Assouad spectrum.
Some of these properties will be used in the proof of \cref{t:moran-exist}, so we cannot formally depend on the corresponding results for the Assouad spectrum.
We draw on ideas from \cite{fhh+2019,fy2018b}.
\begin{proposition}\label{p:a-props}
    Let $\varphi\in\mathcal{A}_d$ be arbitrary.
    Then the following properties hold:
    \begin{enumerate}[nl,r]
        \item\label{im:ep-exist} The limits $\varphi(0)\coloneqq\lim_{\theta\to 0}\varphi(\theta)$ and $\varphi(1)\coloneqq\lim_{\theta\to 1}\varphi(\theta)$ exist.
        \item\label{im:uniform-ct} Each $\varphi\in\mathcal{A}_d$ is uniformly continuous.
        \item\label{im:ep-extreme} $\varphi(0)=\inf_{\theta\in(0,1)}\varphi(\theta)$ and $\varphi(1)=\sup_{\theta\in(0,1)}\varphi(\theta)$.
        \item\label{im:tail-const} For any $\theta_0\in(0,1)$, if $\varphi(\theta_0)=\varphi(1)$, then $\varphi(\theta_0)=\varphi(\theta)$ for all $\theta_0<\theta<1$.
        \item\label{im:zero-bd} If $\varphi(0)=0$, then $\varphi(\theta)=0$ for all $\theta$.
    \end{enumerate}
\end{proposition}
\begin{proof}
    First, we show that $\varphi(\theta)$ is continuous on $(0,1)$.
    For $0<\theta_1<\theta_2<1$ we have $\theta_1<\theta_1/\theta_2<1$, so applying \cref{e:assouad-upper-bd} we obtain
    \begin{equation}\label{e:a-cont-bound}
        (1-\theta_2)\varphi(\theta_2)\leq (1-\theta_1)\varphi(\theta_1)\leq\frac{\theta_1}{\theta_2}(1-\theta_2)\varphi(\theta_2)+\Bigl(1-\frac{\theta_1}{\theta_2}\Bigr)\varphi\Bigl(\frac{\theta_1}{\theta_2}\Bigr).
    \end{equation}
    This implies that
    \begin{equation*}
        \varphi(\theta_1)-\varphi(\theta_2)\leq\frac{\varphi(\theta_1/\theta_2)}{\theta_2(1-\theta_1)}(\theta_2-\theta_1).
    \end{equation*}
    Similarly, from the first inequality of \cref{e:a-cont-bound},
    \begin{equation*}
        \varphi(\theta_2)-\varphi(\theta_1)\leq\left(\frac{1-\theta_1}{1-\theta_2}-1\right)\varphi(\theta)=\frac{\theta_2-\theta_1}{1-\theta_1}\varphi(\theta_1).
    \end{equation*}
    Since $\varphi(\theta_1/\theta_2)\leq d$ and $\varphi(\theta_1)\leq d$, it follows that $\varphi(\theta)$ is Lipschitz on any closed subinterval of $(0,1)$, and therefore continuous on $(0,1)$.

    Now consider \cref{im:ep-exist}.
    Observe that $(1-\theta)\varphi(\theta)$ is a bounded decreasing function of $\theta$, so $\lim_{\theta\to 0}(1-\theta)\varphi(\theta)$ exists so $\lim_{\theta\to 0}\varphi(\theta)$ exists as well.
    To see that $\lim_{\theta\to 1}\varphi(\theta)$ exists, we use the proof from \cite[Section~3.2]{fhh+2019}.
    Set $L=\limsup_{\theta\to 1}\varphi(\theta)$ and let $\varepsilon>0$.
    Since $\varphi(\theta)$ is continuous, we can find $0<u<v<1$ such that $\varphi(\theta)>L-\varepsilon$ for all $\theta\in[u,v]$.
    Thus by \cref{l:nth-root-bd}, with
    \begin{equation*}
        X\coloneqq\bigcup_{n=1}^\infty[u^{1/n},v^{1/n}]
    \end{equation*}
    we have $\varphi(\theta)>L-\varepsilon$ for all $\theta\in X$.
    But $v^{1/n}\geq u^{1/(n+1)}$ for all $n\geq n_0$ with $\frac{n_0}{n_0+1}\geq\frac{\log v}{\log u}$, so in fact $(u^{1/n_0},1)\subset X$.
    Thus $\lim_{\theta\to 1}\varphi(\theta)$ exists as well.
    In particular, combining the existence of endpoint limits with continuity of $\varphi$ on $(0,1)$, \cref{im:uniform-ct} also follows immediately.

    To see \cref{im:ep-extreme}, if $\theta_1\in(0,1)$, then $\theta_n=\theta_1^{1/n}$ is a sequence converging monotonically to $1$ with $\varphi(\theta_n)\geq \varphi(\theta_1)$ by \cref{l:nth-root-bd}.
    Thus $\varphi(1)\geq \varphi(\theta_1)$.
    Similarly $\varphi(\theta_1^n)\leq \varphi(\theta_1)$ for any $n\in\N$, and $\lim_{n\to\infty}\theta_1^n=0$.
    But $\theta_1$ was arbitrary, giving \cref{im:ep-extreme}.

    Now we see \cref{im:tail-const}.
    Suppose $\varphi(1)=\varphi(\theta_1)$ for some $0<\theta_1<1$.
    By \cref{e:a-cont-bound},
    \begin{equation*}
        (1-\theta_1)\varphi(1)-(1-\theta_2)\varphi(\theta_2)\leq(\theta_2-\theta_1)\varphi(\theta_1/\theta_2)\leq(\theta_2-\theta_1)\varphi(1)
    \end{equation*}
    since $\varphi(\theta_1/\theta_2)\leq \varphi(1)$ by \cref{im:ep-extreme}.
    This implies that $\varphi(1)\leq \varphi(\theta_2)$, so \cref{im:tail-const} follows.

    To see \cref{im:zero-bd}, if $\varphi(0)=0$, then $\lim_{\theta\to 0}(1-\theta)\varphi(\theta)=0$.
    But $(1-\theta)\varphi(\theta)$ is a decreasing function of $\theta$, so $(1-\theta)\varphi(\theta)=0$ for all $\theta\in(0,1)$, i.e.~$\varphi(\theta)=0$ for all $\theta\in(0,1)$.
\end{proof}

\subsection{Rate constraints and increasing functions}\label{ss:rate}
Now, we obtain bounds on growth rates of functions in $\mathcal{A}_d$.
We recall that the Dini derivative is defined in \cref{e:Dini-intro}.
We obtain the following regularity property for functions $\varphi\in\mathcal{A}_d$.
\begin{proposition}\label{p:lip-bounds}
    Let $\varphi\in\mathcal{A}_d$ be arbitrary and $\theta\in(0,1)$.
    Then
    \begin{equation*}
        -\frac{\varphi(1)-\varphi(\theta)}{1-\theta}\leq\diniu+\varphi(\theta)\leq\frac{\varphi(\theta)}{1-\theta}
    \end{equation*}
    In particular, $\varphi$ is $d/\delta$-Lipschitz on $[0,1-\delta]$ for any $\delta>0$.
\end{proposition}
\begin{proof}
    The first inequality in \cref{e:assouad-upper-bd} is equivalent to saying that $\beta(\theta)=(1-\theta)\varphi(\theta)$ is decreasing.
    Since $\varphi$ is continuous by \cref{p:a-props}~\cref{im:uniform-ct}, by \cite[Corollary~11.4.2]{bru1994} $\beta$ is decreasing if and only if $\diniu+\beta(\theta)=-\varphi(\theta)+(1-\theta)\diniu+\varphi(\theta)\leq 0$, or equivalently
    \begin{equation*}
        \diniu{+}\varphi(\theta)\leq \frac{\varphi(\theta)}{1-\theta}.
    \end{equation*}
    This gives the upper bound.

    To obtain the lower bound, let $0<\lambda<\theta<1$ be arbitrary.
    By \cref{e:assouad-upper-bd},
    \begin{align*}
        -\varphi(\lambda/\theta)\leq\frac{\beta(\lambda)-\beta(\theta)}{\lambda-\theta}
    \end{align*}
    and taking $\theta\to\lambda$ from the right,
    \begin{equation*}
        -\varphi(1)\leq \diniu+\beta(\lambda)=-\varphi(\theta)+(1-\theta)\diniu+\varphi(\theta).
    \end{equation*}
    Since $0\leq \varphi(\theta)\leq d$ and $0\leq \varphi(1)-\varphi(\theta)\leq d$, it follows that $\varphi$ is $d/\delta$-Lipschitz on $[0,1-\delta]$ for any $\delta>0$.
\end{proof}
\begin{remark}
    In \cref{ss:hol-failure}, we will see that, in general, elements of $\mathcal{A}_d$ need not be Lipschitz (in fact, not even Hölder) on the entire interval $[0,1]$.
\end{remark}
Now, we obtain the following result concerning increasing functions.
\begin{lemma}\label{l:incr-bound}
    If $\varphi\colon(0,1)\to[0,d]$ is increasing, then $\varphi\in\mathcal{A}_d$ if and only if
    \begin{equation}\label{e:dini-bd}
        \diniu+ \varphi(\theta)\leq \frac{\varphi(\theta)}{1-\theta}.
    \end{equation}
\end{lemma}
\begin{proof}
    The forward direction is \cref{p:lip-bounds}.
    To obtain the reverse implication, let $0<\lambda<\theta<1$.
    Since $\varphi$ is increasing, if $\theta\leq\lambda/\theta$, then $\varphi(\lambda)\leq \varphi(\theta)\leq \varphi(\lambda/\theta)$ and
    \begin{equation*}
        0\leq(1-\lambda)\varphi(\lambda)-(1-\theta)\varphi(\theta)\leq(\theta-\lambda)\varphi(\theta)\leq(\theta-\lambda)\varphi\Bigl(\frac{\lambda}{\theta}\Bigr)
    \end{equation*}
    and if $\lambda/\theta\leq\theta$,
    \begin{equation*}
        0\leq(1-\lambda)\varphi(\lambda)-(1-\theta)\varphi(\theta)\leq(1-\lambda)\varphi\Bigl(\frac{\lambda}{\theta}\Bigr)-(1-\theta)\varphi\Bigl(\frac{\lambda}{\theta}\Bigr)=(\theta-\lambda)\varphi\Bigl(\frac{\lambda}{\theta}\Bigr)
    \end{equation*}
    which is \cref{e:assouad-upper-bd}.
\end{proof}
We obtain the following convenient application, which we use to characterize the upper Assouad spectra.
\begin{corollary}\label{c:upper-a}
    Let $\varphi\in\mathcal{A}_d$.
    Then $\overline{\varphi}\in\mathcal{A}_d$ where
    \begin{equation*}
        \overline{\varphi}(\theta)=\sup_{0<\theta'\leq\theta}\varphi(\theta').
    \end{equation*}
\end{corollary}
\begin{proof}
    As proven in \cref{l:incr-bound}, since $\overline{\varphi}(\theta)$ is increasing, we only need to verify that $\diniu+\overline{\varphi}(\theta)\leq\overline{\varphi}(\theta)/(1-\theta)$.
    Since $\varphi(\theta)\leq\overline{\varphi}(\theta)$, it suffices to show $\diniu+\overline{\varphi}\leq\max\{\diniu+ \varphi,0\}$.

    Fix $\theta_0$ and let $(\theta_n)_{n=1}^\infty\to\theta_0$ be strictly decreasing.
    Passing to a subsequence if necessary, we may assume $\overline{\varphi}(\theta_n)>\overline{\varphi}(\theta_0)$ for all $n$; otherwise $\diniu+\overline{\varphi}(\theta_0)\leq 0$.
    Thus for each $n$ there is $\theta_0<\theta_n'\leq\theta_n$ be such that $\varphi(\theta_n')=\overline{\varphi}(\theta_n)$.
    Thus
    \begin{equation*}
        \frac{\overline{\varphi}(\theta_n)-\overline{\varphi}(\theta_0)}{\theta_n-\theta_0}\leq\frac{\varphi(\theta'_n)-\varphi(\theta_0)}{\theta_n'-\theta_0}\leq\diniu+\varphi(\theta_0).
    \end{equation*}
    But $(\theta_n)_{n=1}^\infty$ was an arbitrary sequence, so the result follows.
\end{proof}

\section{Classifying the forms of Assouad spectra}\label{s:class}
In this section, we prove the main classification result, \cref{it:assouad-char}.
\subsection{Bounding the Assouad spectrum}
We recall the following general bounds, which are given in \cite[Proposition~3.4]{fy2018b} and \cite[Theorem~3.3.1]{fra2021}.
We include the details here for completeness.
\begin{proposition}\label{p:upper}
    For any set $F\subset\R^d$, the function $\varphi(\theta)=\dimAs\theta F$ is in $\mathcal{A}_d$.
\end{proposition}
\begin{proof}
    Let $0<\theta_1<\theta_2<1$ and let $\varepsilon>0$ be arbitrary.
    For $\delta>0$ sufficiently small, since $B(x,\delta^{\theta_2})\subset B(x,\delta^{\theta_1})$ for all $x\in F$,
    \begin{align*}
        \sup_{x\in F}N_\delta(F\cap B(x,\delta^{\theta_1}))&\geq\sup_{x\in F}N_\delta(F\cap B(x,\delta^{\theta_2}))\\
                                                           &\geq\left(\frac{\delta^{\theta_2}}{\delta}\right)^{(\varphi(\theta_2)-\varepsilon)}\\
                                                           &=\left(\delta^{\theta_1-1}\right)^{(\varphi(\theta_2)-\varepsilon)\left(\frac{1-\theta_2}{1-\theta_1}\right)}
    \end{align*}
    which proves that $(1-\theta_1)\varphi(\theta_1)\geq(1-\theta_2)(\varphi(\theta_2)-\varepsilon)$.
    This gives the lower inequality in \cref{e:assouad-upper-bd}.

    To obtain the upper inequality, by covering $B(x,\delta^{\theta_1})$ by balls with radius $\delta^{\theta_2}$,
    \begin{equation*}
        \sup_{x\in F} N_\delta(F\cap B(x,\delta^{\theta_1}))\leq\sup_{x\in F}N_{\delta^{\theta_2}}(F\cap B(x,\delta^{\theta_1}))\sup_{x\in F}N_\delta(F\cap B(x,\delta^{\theta_2})).
    \end{equation*}
    This implies for all $\delta>0$ sufficiently small
    \begin{align*}
        \sup_{x\in F}N_\delta(F\cap B(x,\delta^{\theta_1}))&\leq\left(\frac{\delta^{\theta_1}}{\delta^{\theta_2}}\right)^{\varphi(\theta_1/\theta_2)+\varepsilon}\left(\frac{\delta^{\theta_2}}{\delta}\right)^{\varphi(\theta_2)+\varepsilon}\\
                                                           &=\left(\delta^{\theta_1-1}\right)^{(\varphi(\theta_1/\theta_2)+\varepsilon)\left(\frac{\theta_2-\theta_1}{1-\theta_1}\right)+(\varphi(\theta_2)+\varepsilon)\left(\frac{1-\theta_2}{1-\theta_1}\right)}
    \end{align*}
    which implies that
    \begin{equation*}
        (1-\theta_1)\varphi(\theta_1)\leq(\theta_2-\theta_1)(\varphi(\theta_1/\theta_2)+\varepsilon)+(1-\theta_2)(\varphi(\theta_2)+\varepsilon)
    \end{equation*}
    as required.
\end{proof}

\subsection{Constructing sets with prescribed spectra}
Now for any $\varphi\in\mathcal{A}_d$ we construct a homogeneous Moran set $C$ such that $\dimAs\theta C=\varphi(\theta)$ for all $\theta\in(0,1)$.
The techniques here are based on ideas first introduced by the author and Banaji used to solve an analogous question for the \emph{intermediate dimensions} \cite{zbl:1448.28009}.
We refer the reader to the paper \cite{zbl:1509.28005} for more details on this general technique.

We first recall the notion of homogeneous Moran sets from \cite{zbl:1509.28005}.
The construction is analogous to the usual $2^d$-corner Cantor set, except that the subdivision ratios need not be the same at each level.

Let $\mathcal{I}=\{0,1\}^d$, set $\mathcal{I}^*=\bigcup_{n=0}^\infty\mathcal{I}^n$, and denote the word of length $0$ by $\varnothing$.
Let $\bm{r}=(r_n)_{n=1}^\infty\subset(0,1/2]$ and for each $n$ and $\bm{i}\in\mathcal{I}$, define $S^n_{\bm{i}}\colon\R^d\to\R^d$ by
\begin{equation*}
    S^n_{\bm{i}}(x)\coloneqq r_n x+b^n_{\bm{i}}
\end{equation*}
where $b^n_{\bm{i}}\in\R^d$ has
\begin{equation*}
    (b^n_{\bm{i}})^{(j)} =
    \begin{cases}
        0 &: \bm{i}^{(j)}=0\\
        1-r_n &: \bm{i}^{(j)}=1
    \end{cases}
    .
\end{equation*}
Given $\sigma=(\bm{i}_1,\ldots,\bm{i}_n)\in\mathcal{I}^n$, we write $S_\sigma=S^1_{\bm{i}_1}\circ\cdots\circ S^n_{\bm{i}_n}$.
Then set
\begin{equation*}
    C=C(\bm{r})\coloneqq\bigcap_{n=1}^\infty \bigcup_{\sigma\in\mathcal{I}^n}S_\sigma([0,1]^d).
\end{equation*}
We refer to the set $C$ as a \defn{homogeneous Moran set}.

Given $\delta>0$, let $k=k(\delta)$ be such that $r_1\cdots r_k\leq \delta < r_1\cdots r_{k-1}$.
We then define
\begin{equation*}
    s(\delta)=s_{\bm{r}}(\delta)\coloneqq\frac{k(\delta)\cdot d\log 2}{-\log\delta}.
\end{equation*}
Heuristically, $s(\delta)$ is the best candidate for the box dimension of $C$ at scale $\delta$.

We now define a family of functions which we may interpret as a reparametrization of the space of sequences $(0,1/2]^{\N}$.
That this forms an alternative representation is described precisely in \cite[Lemma~3.4]{zbl:1509.28005}.
\begin{definition}\label{d:g}
    Let $0\leq \lambda\leq\alpha\leq d$ and let $\mathcal{G}(\lambda,\alpha)$ denote the set of functions $g\colon\R\to(\lambda,\alpha)$ satisfying
    \begin{equation*}
        \lambda-(\lambda-g(y))\exp(-t)\leq g(y+t)\leq \alpha-(\alpha-g(y))\exp(-t)
    \end{equation*}
    for any $y\in\R$ and $t>0$.
\end{definition}
This family is the same as the family defined in \cite[Definition~3.1]{zbl:1509.28005} (the argument is given in \cite[Lemma~3.2]{zbl:1509.28005}).

To construct sets with prescribed Assouad spectrum, we will follow the approach from \cite{zbl:1509.28005} and use homogeneous Moran sets.
Given a function $g\in\mathcal{G}(\lambda,\alpha)$ and $w\in\R$, we define the \emph{offset} $\kappa_w(g)\in\mathcal{G}(\lambda,\alpha)$ by
\begin{equation*}
    \kappa_w(g)(x)=\begin{cases}
        g(x-w) &: x\geq w,\\
        g(0) &: x\leq w.
    \end{cases}
\end{equation*}
We also say that a function $g\in\mathcal{G}(\lambda,\alpha)$ is \emph{rapidly decreasing} if there is a $y\in\R$ and a constant $C>0$ so that for all $x\geq y$,
\begin{equation}\label{e:rapidly-decreasing}
    g(x)\leq g(y)\exp(y-x)+C\exp(-x).
\end{equation}
Note that if $g$ is rapidly decreasing, then $\lim_{x\to\infty}g(x)=0$.
Moreover, for all $w\in\R$, $g$ is not rapidly decreasing if and only if $\kappa_w(g)$ is not rapidly decreasing.

The following lemma is slightly different than \cite[Lemma~3.4]{zbl:1509.28005} but follows by essentially the same proof.
\begin{lemma}\label{l:moran}
    Let $0\leq\lambda\leq\alpha\leq d$ and let $\tilde{g}\in\mathcal{G}(\lambda,\alpha)$.
    Suppose $\tilde{g}$ is not rapidly decreasing.
    Then there is a constant $w_0\in\R$ depending only on $\tilde{g}(0)$ and $d$ such that for all $w\geq w_0$, there exists a sequence $\bm{r}\coloneqq (r_j)_{j=1}^\infty\subset(0,1/2]$ so that $g\coloneqq\kappa_{w}(\tilde{g})$ satisfies
    \begin{equation}\label{e:g-disc}
        |s_{\bm{r}}(\exp(-\exp(x)))-g(x)|\leq d\log(2)\cdot\exp(-x)
    \end{equation}
    for all $x\geq w_0$.
\end{lemma}
Using \cref{l:moran}, we establish the following general result which allows us to prescribe Assouad spectra for homogeneous Moran sets.
\begin{proposition}\label{p:Moran-cnst}
    Let $d\in\N$ and $g\in\mathcal{G}(0,d)$.
    Then there exists a homogeneous Moran set $C$ such that
    \begin{equation}\label{e:g-formula}
        \dimAs\theta C=\limsup_{x\to \infty}\frac{g\Bigl(x+\log\frac{1}{\theta}\Bigr)-\theta g(x)}{1-\theta}.
    \end{equation}
\end{proposition}
\begin{proof}
    If $g$ is rapidly decreasing, then $\lim_{x\to\infty}g(x)=0$.
    By a short argument, this implies that
    \begin{equation*}
        \limsup_{x\to\infty} \frac{g\Bigl(x+\log\frac{1}{\theta}\Bigr) - \theta g(x)}{1-\theta} = 0.
    \end{equation*}
    Thus we can define the Moran set $C(\bm{r})$ where $\bm{r}$ is a sequence converging to 0.

    Otherwise, performing an appropriate translation of $g$ which does not change \cref{e:g-formula}, \cref{l:moran} provides a sequence $\bm{r}\subset(0,1/2]$ satisfying \cref{e:g-disc}.
    To obtain \cref{e:g-formula}, by definition of the Assouad spectrum,
    \begin{equation*}
        \dimAs\theta C=\limsup_{\delta\to 0}\sup_{x\in C}\frac{\log N_{\delta^{1/\theta}}(C\cap B(x,\delta))}{(1-1/\theta)\log\delta}.
    \end{equation*}
    Observe that there is some constant $M>0$ such that $B(x,\delta)$ intersects at most $M$ cylinders in level $k(\delta)$.
    In particular, $C\cap B(x,\delta)$ can be covered by $M\cdot 2^{d(k(\delta^{1/\theta})-k(\delta))}$ balls of radius $\delta^{1/\theta}$.
    On the other hand, $C\cap B(x,\delta)$ contains an interval in level $k(\delta)$, and therefore contains a $\delta$-separated subset of size $2^{d(k(\delta^{1/\theta})-1-k(\delta))}$.
    Thus there is a constant $M'>0$ so that
    \begin{equation*}
        M'\cdot 2^{d(k(\delta^{1/\theta})-k(\delta))}\leq \sup_{x\in C}N_{\delta^{1/\theta}}(C\cap B(x,\delta))\leq M\cdot 2^{d(k(\delta^{1/\theta})-k(\delta))}
    \end{equation*}
    and therefore
    \begin{align*}
        \limsup_{\delta\to 0}\sup_{x\in C}\frac{\log N_{\delta^{1/\theta}}(C\cap B(x,\delta))}{(1-1/\theta)\log\delta} &=\limsup_{\delta\to 0}\frac{\theta(k(\delta^{1/\theta})-k(\delta))\cdot d\log 2}{(1-\theta)\cdot(-\log\delta)}\\
                                                                                                                       &=\limsup_{\delta\to 0}\frac{s(\delta^{1/\theta})-\theta \cdot s(\delta)}{1-\theta}.
    \end{align*}
    Taking $\delta>0$ small and applying \cref{e:g-disc} yields the desired formula.
\end{proof}
\begin{definition}
    Given a sequence of continuous functions $(f_k)_{k=1}^\infty$ each defined on some interval $[0,a_k]$, the \defn{concatenation} of $(f_k)_{k=1}^\infty$ is the function
    \begin{equation*}
        f\colon(-\infty,\sum_{k=1}^\infty a_k)\to\R
    \end{equation*}
    given as follows: for each $x>0$ with $\sum_{j=0}^{k-1} a_j<x\leq\sum_{j=0}^{k}a_j$ where $a_0=0$ we define
    \begin{equation*}
        f(x)=f_k\left(x-\sum_{j=0}^{k-1} a_j\right),
    \end{equation*}
    and we define $f(x)=f_1(0)$ for $x\leq 0$.
    The concatenation of a finite tuple of functions is defined similarly.
\end{definition}
We next prove the converse direction of \cref{it:assouad-char}.
For the convenience of the reader, we also give an explicit description of the construction technique in $\R$.
Note that in the proof of \cref{t:moran-exist}, the precise choice of the contractions $(r_i)_{i=1}^\infty\subset(0,1/2]$ is concealed in the application of \cref{l:moran} in \cref{p:Moran-cnst}.

Let $\varphi\in\mathcal{A}_1$ be some fixed function.
Fix some small constant $\delta_1$.
Then we will inductively choose constants $r^{(n)}_1,\ldots,r^{(n)}_{m_n}$ in $(0,1/2]$ for each $n\in\N$ so that for each $1\leq j\leq m_n$
\begin{equation}\label{e:rn-choice}
    2^j\approx\Bigl(\frac{1}{r^{(n)}_1\cdots r^{(n)}_j}\Bigr)^{\varphi(\theta)}
\end{equation}
where $\theta$ is such that $\delta_n^{1/\theta}\approx\delta_n\cdot r^{(n)}_1\cdots r^{(n)}_j$, and $m_n$ satisfies $\delta_n\cdot r^{(n)}_1\cdots r^{(n)}_{m_n}\approx \delta_n^{n}$.
Then take $R_n$ very small, and set $\delta_{n+1}=\delta_n\cdot r^{(n)}_1\cdots r^{(n)}_{m_n}\cdot R_n$.
Now let $C$ denote the Moran set corresponding to the sequence
\begin{equation*}
    (\delta_1,r^{(1)}_1,\ldots,r^{(1)}_{m_1},R_1,r^{(2)}_1,\ldots,r^{(2)}_{m_2},R_2,\ldots).
\end{equation*}
For $n\in\N$ and $x\in C$, $N_{\delta_n^{1/\theta}}(C\cap B(x,\delta_n))\approx 2^j$ where $\delta_n^{1/\theta}\approx\delta_n\cdot r^{(n)}_1\cdots r^{(n)}_j$.
In particular, \cref{e:rn-choice} guarantees that the Assouad spectrum of $C$ with respect to $\theta$ at scale $\delta_n$ is precisely $\varphi(\theta)$, for all $n$ sufficiently large so that $1/n\leq\theta$.

The main details of the proof are to show (1) that such a choice of the constants $r_i$ is possible, and (2) that for fixed $\theta$ and sufficiently small scales $\delta$ not of the form $\delta_n$, the Assouad spectrum at $\theta$ of $C$ at scale $\delta$ is at most $\varphi(\theta)$.
\begin{theorem}\label{t:moran-exist}
    Let $\varphi\in\mathcal{A}_d$ be arbitrary.
    Let $\alpha$ be such that $\varphi(1)\leq\alpha\leq d$.
    Then there exists a homogeneous Moran set $C\subseteq\R^d$ such that $\dimA C=\alpha$ and, for all $\theta\in(0,1)$,
    \begin{equation*}
        \dimAs\theta C=\varphi(\theta).
    \end{equation*}
\end{theorem}
\begin{proof}
    We may assume $\alpha>0$, or the result is immediate.
    We will prove the result for the Assouad spectrum, and then explain how to modify the proof to accommodate the Assouad dimension as well.

    First, we apply some convenient rescaling to $\varphi(\theta)$.
    Given $y\in (0,\infty)$, $\exp(-y)\in(0,1)$ so we may define
    \begin{equation*}
        \xi(y)=(1-\exp(-y))\varphi(\exp(-y)).
    \end{equation*}
    In particular, given $0<y_1<y_2<\infty$, it follows that $0<\exp(-y_2)<\exp(-y_1)<1$ so
    \begin{align*}
        0&\leq (1-\exp(-y_2))\varphi(\exp(-y_2))-(1-\exp(-y_1))\varphi(\exp(-y_1))\\
         &\leq \exp(-y_1)\left(1-\exp(-(y_2-y_1))\right)\varphi(\exp(-(y_2+y_1)))
    \end{align*}
    or equivalently
    \begin{equation}\label{e:phi-props}
        0\leq\xi(y_2)-\xi(y_1)\leq\exp(-y_1)\xi(y_2-y_1).
    \end{equation}
    Moreover, observe that $\varphi(1)=\lim_{y\to 0}\varphi(\exp(-y))$ so $\lim_{y\to 0}\xi(y)=0$, and similarly $\lim_{y\to\infty}\xi(y)=\varphi(0)$.
    In particular $\xi$ is continuous, increasing, and bounded.

    Now for $z\in(0,\alpha)$, let $\xi_z$ denote the function
    \begin{equation*}
        \xi_z(y)=\xi(y)+\exp(-y)z
    \end{equation*}
    and similarly $\Psi_z(y)=\exp(-y)z$.
    We note that $\xi_z(0)=\Psi_{z}(0)=z$.

    Now, let $z_1=0$ and choose constants $w_n,z_n$ such that the functions $f_n\coloneqq \xi_{z_n}|_{[0,n]}$ and $e_n\coloneqq\Psi_{w_n}|_{[0,n]}$ satisfy $f_n(n)=e_{n}(0)$ and $e_n(n)=f_{n+1}(0)$ for all $n\in\N$.
    Then, let $g$ be the infinite concatenation of the sequence
    \begin{equation*}
        \bigl(f_1,e_1,f_2,e_2,\ldots\bigr).
    \end{equation*}
    This construction is illustrated in \cref{f:concat-const}.
    \begin{figure}[t]
        \begin{center}
            \begin{tikzpicture}[>=stealth,scale=3.6,font=\tiny]
                \draw[thick,->] (-0.05,0) -- (4.1,0);
                \draw[thick,->] (0,-0.05) -- (0,0.62);

                \begin{scope}[thick, smooth, variable=\y]
                    \draw[domain=0:ln(3/2)] plot ({\y}, {1-exp(-\y)});
                    \draw[domain=ln(3/2):ln(2)] plot ({\y}, {1/3});
                    \draw[domain=ln(2):1] plot ({\y}, {5/6-exp(-\y)});

                    \draw[domain=1:2, gray] plot ({\y}, {(5*exp(1)/6-1)*exp(-\y)});

                    \draw[domain=0:ln(3/2)] plot ({\y+2}, {(5/(6*exp(1))-1/(exp(2)))*exp(-\y)+1-exp(-\y)});
                    \draw[domain=ln(3/2):ln(2)] plot ({\y+2}, {(5/(6*exp(1))-1/(exp(2)))*exp(-\y)+1/3});
                    \draw[domain=ln(2):ln(3)] plot ({\y+2}, {(5/(6*exp(1))-1/(exp(2)))*exp(-\y)+5/6-exp(-\y)});
                    \draw[domain=ln(3):2] plot ({\y+2}, {(5/(6*exp(1))-1/(exp(2)))*exp(-\y)+1/2});
                \end{scope}

                \draw[dotted] (1,-0.11) -- (1,0.46545389);
                \draw[dotted] (2,-0.11) -- (2,0.17123092);
                \draw[dotted] (4,-0.11) -- (4,0.52317358);

                \draw[<->] (0+0.01,-0.1) -- node[fill=white]{$f_1$} (1-0.01,-0.1);
                \draw[<->] (1+0.01,-0.1) -- node[fill=white]{$e_1$} (2-0.01,-0.1);
                \draw[<->] (2+0.01,-0.1) -- node[fill=white]{$f_2$} (4-0.01,-0.1);
            \end{tikzpicture}
        \end{center}
        \caption{The concatenation of $(f_1,e_1,f_2)$ corresponding to a function $\phi\in\mathcal{C}_d$ defined in \cref{ss:non-mono-ex} restricted to the domain $(0,\infty)$.}
        \label{f:concat-const}
    \end{figure}

    First, let us verify that $g\in\mathcal{G}(0,\varphi(1))\subset\mathcal{G}(0,\alpha)$.
    Since membership of $\mathcal{G}$ is equivalent to a pointwise derivative constraint (see \cite[Lemma~3.2]{zbl:1509.28005}), it suffices to verify \cref{d:g} piecewise.
    Let $n\in\N$.
    Note that $e_n\in\mathcal{G}(0,\varphi(1))$ since the $e_n$ are differentiable with $e_n'(x)=\varphi(0)-e_n(x)$.
    Next let $0<y<y+t<\infty$.
    First observe that
    \begin{equation*}
        \xi(y+t)\leq\xi(t)+\xi(y)\exp(-t)\leq(1-\exp(-t))\varphi(1)+\xi(y)\exp(-t)
    \end{equation*}
    by \cref{e:phi-props} and \cref{p:a-props}~\cref{im:ep-extreme}.
    Thus
    \begin{align*}
        f_n(y+t)&=\xi(y+t)+\exp(-(y+t))z_n\\
                &\leq (1-\exp(-t))\varphi(1)+\xi(y)\exp(-t) +\exp(-(y+t))z_n\\
                &= (1-\exp(-t))\varphi(1)+f_n(y)\exp(-t)
    \end{align*}
    as required.
    To obtain the other bound, since $\xi$ is increasing,
    \begin{align*}
        f_n(y+t) &= \xi(y+t)+\exp(-(y+t))z_n\\
                 &\geq \xi(y)\exp(-t)+\exp(-y)\exp(-t)z_n\\
                 &= f_n(y)\exp(-t).
    \end{align*}

    Now, let $C$ denote the Moran set corresponding to the function $g$.
    Let $\theta\in(0,1)$: we must show that $\dimAs\theta C=\varphi(\theta)$.
    Let $\tau=\log(1/\theta)$.
    By \cref{p:Moran-cnst}, it suffices to show
    \begin{equation}\label{e:limsup-form}
        \varphi(\theta)=\limsup_{x\to\infty}\frac{g(x+\tau)-\theta g(x)}{1-\theta}.
    \end{equation}
    For $n\in\N$ set $x_n=2\sum_{i=1}^{n-1}i$ and let $N\in\N$ be sufficiently large so that $N\geq \tau+1$.

    Now if $n\geq N$, $g(x_n+\tau)=f_n(\tau)$ and $g(x_n)=f_n(0)=z_n$ so that
    \begin{equation*}
        \frac{g(x_n+\tau)-\theta g(x_n)}{1-\theta} = \frac{(1-\theta)\varphi(\theta)+\theta z_n-\theta z_n}{1-\theta}=\varphi(\theta).
    \end{equation*}
    This gives the lower bound in \cref{e:limsup-form}.

    It remains to see the upper bound.
    We first observe for all $y>0$ and $z\in\R$ that
    \begin{equation*}
        \frac{\xi_z(y+\tau)-\theta \xi_z(y)}{1-\theta}\leq \varphi(\theta).
    \end{equation*}
    Indeed, expanding the definition of $\xi_z$ and applying \cref{e:phi-props},
    \begin{align*}
        \xi_z(y+\tau)-\theta \xi_z(y) &=\xi(y+\tau)+\exp(-(y+\tau))z-\exp(-\tau)(\xi(y)+\exp(-y)z)\\
                               &= \xi(y+\tau)-\exp(-\tau)\xi(y)\\
                               &\leq\xi(\tau)=(1-\theta)\varphi(\theta).
    \end{align*}

    Now let $x\geq x_{N}$ be arbitrary and let $n$ be such that $x\in[x_n-(n-1),x_n+n]$.
    First note that for $y\in[x_n-(n-1),x_n+2n]$, $g(y)=\exp(-(y-x_n))z_n+\phi(y)$ where
    \begin{equation*}
        \phi(y)=
        \begin{cases}
            0 &: x_n-(n-1)\leq y\leq x_n\\
            \xi(y-x_n) &: x_n\leq y\leq x_n+n\\
            \xi(n)\exp(-(y-x_n+n)) &: x_n+n\leq y\leq x_n+2n
        \end{cases}
    \end{equation*}
    by choice of the constants $w_n$ and $z_n$.
    If $x\in[x_n,x_n+n]$, since $x+\log(1/\theta)\leq x_n+2n$ and $g(y)\leq\xi_{z_n}(y-x_n)$ for all $y\in[x_n,x_n+2n]$, the prior computation shows that $g(x+\tau)-\theta g(x)\leq(1-\theta)\varphi(\theta)$.
    Otherwise, $x\in[x_n-(n-1),x_n]$.
    If $x+\tau\leq x_n$, then $g(x+\tau)-\theta g(x) = 0\leq(1-\theta)\varphi(\theta)$, and if $x_n<x+\tau\leq x_n+n$, then
    \begin{align*}
        g(x+\tau)-\theta g(x)=\xi(x+\tau-x_n)\leq\xi(\tau)
    \end{align*}
    since $\xi$ is increasing.
    Thus \cref{e:limsup-form} holds, finishing the proof.

    In order to obtain the result for the Assouad dimension as well, we modify the construction as follows.
    Define functions $u_n\colon[0,1/n]\to(0,\alpha)$ by the rule $u_n(x)=\alpha-(\alpha-q_n)\exp(-x)$.
    Choosing the constants $q_n$ appropriately and modifying the constants $w_n$ and $z_n$, the concatenation $\tilde{g}$ of the sequence
    \begin{equation*}
        (f_1,e_1,u_1,f_2,e_2,u_2,\ldots)
    \end{equation*}
    is continuous and $\tilde{g}\in G(0,\alpha)$ since $\alpha\geq \varphi(1)$.
    Since the $u_n$ are supported on intervals with lengths converging to $0$, the same arguments as before yield the correct bounds for $\dimAs\theta C$ up to an error decaying to 0 as $n$ goes to infinity.
    On the other hand, the same arguments as given in \cite[Lemma~3.7~and~Theorem~3.9]{zbl:1509.28005} give that $\dimA C=\alpha$.
    We leave the precise details to the reader.
\end{proof}
\subsection{Families of monotonic and non-monotonic spectra}\label{ss:non-mono-ex}
In this section, we define two general parametrized families of functions in $\mathcal{A}_d$.
The first is a 2-parameter family composed of increasing functions, and the second is a 3-parameter family composed of functions which (outside of degenerate cases) are non-monotonic.
\subsubsection{Monotonic spectra}
Let
\begin{equation*}
    M_d=\bigl\{(\kappa,c):0\leq\kappa\leq d,0<c<1\bigr\}
\end{equation*}
and for $\bm{i}=(\kappa,c)\in M_d$, we may define
\begin{equation*}
    f_{\bm{i}}(\theta) = \begin{cases}
        \kappa(1-c) &: \theta\in[0,c]\\
        \kappa(1-\theta) &: \theta\in[c,1]
    \end{cases}
    .
\end{equation*}
Then, let
\begin{equation}\label{e:Id}
    \mathcal{M}_d\coloneqq\Bigl\{\theta\mapsto \frac{f_{\bm{i}}(\theta)}{1-\theta}:\bm{i}\in M_d\Bigr\}.
\end{equation}
A direct argument shows that each $\varphi\in\mathcal{M}_d$ is an increasing element of $\mathcal{A}_d$.
\subsubsection{Non-monotonic spectra}
This family generalizes the example considered in \cite[Theorem~3.4.16]{fra2021}.
Let
\begin{equation*}
    C_d=\bigl\{(\kappa,c_1,c_2):0\leq \kappa\leq d,0\leq c_1\leq c_2\leq c_1^{1/2}<1\bigr\}.
\end{equation*}
Suppose $\bm{c}=(\kappa,c_1,c_2)\in C_d$.
If $c_1=0$, let $h_{\bm{c}}(\theta)=\kappa(1-\theta)$ for $\theta\in[0,1]$.
Otherwise, $c_2\leq c_1/c_2$.
Thus we may define $h=h_{\bm{c}}\colon[0,1]\to[0,d]$ to be the unique continuous function which has slope $0$ on $[0,c_1]\cup[c_2,c_1/c_2]$, has slope $-\kappa$ on $[c_1,c_2]\cup[c_1/c_2,1]$, and satisfies $h(1)=0$.
Now, let
\begin{equation}\label{e:Fd}
    \mathcal{C}_d=\Bigl\{\theta\mapsto\frac{h_{\bm{c}}(\theta)}{1-\theta}:\bm{c}\in C_d\Bigr\}
\end{equation}
We note that $h_{\bm{c}}$ satisfies a certain rescaling invariance: for $c_2^2\leq\theta\leq c_1$,
\begin{equation*}
    h_{\bm{c}}(\theta)-h_{\bm{c}}\Bigl(\frac{\theta}{c_2}\Bigr)=\kappa(c_2-c_1).
\end{equation*}
In particular, $h_{\bm{c}}(c_2^2)/(1-c_2^2)=h_{\bm{c}}(c_2)/(1-c_2)$.

There are degenerate cases: if $c_2=c_1^{1/2}$, then $h_{(\kappa,c_1,c_2)}=f_{(\kappa,c_1)}$ and if $c_1=c_2$ or $\kappa=0$, then $h_{\kappa,c_1,c_2}=0$.
Otherwise, $h_{\bm{c}}(\theta)/(1-\theta)$ is strictly increasing on $[0,c_1]$ and $[c_2,c_1/c_2]$, constant on $[c_1/c_2,1]$, and strictly decreasing on $[c_1,c_2]$.
A plot of the function $h_{\bm{c}}(\theta)/(1-\theta)$ for non-degenerate parameters is given in \cref{f:decr-plot}.
\begin{figure}[t]
    \begin{center}
        \begin{tikzpicture}[scale=8,font=\tiny,>=stealth]
            \path (0,1/3) -- (1.4,1.05);

            \draw[thick,->] (-0.05,1/3) -- (1.05,1/3);
            \draw[thick,->] (0,1/3-0.05) -- (0,1.05);

            \begin{scope}[thick, dashed, smooth, variable=\x]
                \draw[domain=1/3:1/2] plot ({\x}, {(1/2)/(1-\x)});
                \draw[domain=1/2:2/3] plot ({\x}, {1});
            \end{scope}

            \begin{scope}[thick, smooth, variable=\x]
                \draw[domain=0:1/3] plot ({\x}, {(1/2)/(1-\x)});
                \draw[domain=1/3:1/2] plot ({\x}, {(5/6-\x)/(1-\x)});
                \draw[domain=1/2:2/3] plot ({\x}, {(1/3)/(1-\x)});
                \draw[domain=2/3:1] plot ({\x}, {1});
            \end{scope}

            \matrix [draw,below left,scale=0.8,column sep=0.3em, inner sep=0.4em, align=left] at (current bounding box.north east) {
                \draw[thick, dashed] (0,-0.3) -- (0.71,-0.3); & \node[inner sep=0pt] {$\frac{f_{(\kappa,c_2)}(\theta)}{1-\theta}$}; \\
                \draw[thick] (0,-0.3) -- (0.71,-0.3); & \node[inner sep=0pt] {$\frac{h_{(\kappa,c_1,c_2)}(\theta)}{1-\theta}$}; \\
            };

            \begin{scope}[dotted]
                \draw (1/3,1/3) -- (1/3,3/4);
                \draw (1/4,1/3) node[below]{$c_2^2$}-- (1/4,2/3);
                \draw (1/3,1/3) node[below]{$c_1$}-- (1/3,3/4);
                \draw (1/2,1/3) node[below]{$c_2$}-- (1/2,1);
                \draw (2/3,1/3) node[below]{$\frac{c_1}{c_2}$}-- (2/3,1);
                \draw (1,1/3) node[below]{$1$}-- (1,1);

                \draw (0,1/2) -- (1,1/2) node[right]{$\kappa\left(1-c_1+c_2-\frac{c_1}{c_2}\right)$};
                \draw (1/4,2/3) -- (1,2/3) node[right]{$\frac{\kappa\left(1-\frac{c_1}{c_2}\right)}{1-c_2}$};
                \draw (1/3,3/4) -- (1,3/4) node[right]{$\frac{\kappa\left(1-\frac{c_1}{c_2}+c_2-c_1\right)}{1-c_1}$};
                \node[right] at (1,1) {$\kappa$};
            \end{scope}
        \end{tikzpicture}
    \end{center}
    \caption{A plot of $f_{\bm{i}}(\theta)/(1-\theta)$ and $h_{\bm{c}}(\theta)/(1-\theta)$ where $\bm{i}=(\kappa,c_2)$ and $\bm{c}=(\kappa,c_1,c_2)$.}
    \label{f:decr-plot}
\end{figure}
\begin{proposition}\label{p:bmFd}
    For any $d\in\N$, $\mathcal{C}_d\subset\mathcal{A}_d$.
\end{proposition}
\begin{proof}
    Fix $0<\theta_1,\theta_2<1$ and $\bm{c}=(\kappa,c_1,c_2)\in C_d$.
    We may assume $0<c_1<1$.
    Since $h_{\bm{c}}$ is decreasing, it suffices to show that
    \begin{equation}\label{e:h-bound}
        \frac{h_{\bm{c}}(\lambda)-h_{\bm{c}}(\theta)}{\theta}\leq h_{\bm{c}}\Bigl(\frac{\lambda}{\theta}\Bigr)
    \end{equation}
    for all $0<\lambda<\theta<1$.
    Since \cref{e:h-bound} is invariant under scaling by a positive factor, we may assume $\kappa=1$.
    We prove this result in cases depending on the positions of $\lambda$ and $\theta$.

    If $\lambda\in[0,c_2^2]\cup[c_2,1]$, then $h_{\bm{c}}(\lambda)/(1-\lambda)\leq h_{\bm{c}}(\theta)/(1-\theta)$ for all $\lambda<\theta$.
    In particular, as argued in \cref{l:incr-bound}, \cref{e:h-bound} holds for all such $\lambda$.
    Moreover, suppose that \cref{e:h-bound} holds for the choice $\lambda=c_1$ and all $\lambda<\theta$.
    Since $h_{\bm{c}}$ is the constant function on $[c_2^2,c_1]$, this implies the bound on $[c_2^2,c_1$].
    Moreover, for $\lambda\in[c_1,c_2]$, since $h_{\bm{c}}(c_1/\theta)-h_{\bm{c}}(\lambda/\theta)\leq(c_1-\lambda)/\theta$,
    \begin{equation*}
        h_{\bm{c}}(\lambda)-h_{\bm{c}}(\theta)\leq\theta h_{\bm{c}}\Bigl(\frac{c_1}{\theta}\Bigr)-(c_1-\lambda)\leq\theta h_{\bm{c}}\Bigl(\frac{\lambda}{\theta}\Bigr).
    \end{equation*}
    Thus it suffices to establish \cref{e:h-bound} for $\lambda=c_1$ and $\theta>\lambda$.
    Write $g(\theta)=(h_{\bm{c}}(c_1)-h_{\bm{c}}(\theta))/\theta$: we must show that $g(\theta)\leq h_{\bm{c}}(c_1/\theta)$.
    \begin{enumerate}
        \item If $\theta\in[c_1,c_2]$, then $c_1/\theta\geq c_1/c_2$ and
            \begin{equation*}
                g(\theta)=1-\frac{c_1}{\theta}=h_{\bm{c}}\Bigl(\frac{c_1}{\theta}\Bigr).
            \end{equation*}
        \item If $\theta\in [c_2,c_1/c_2]$ then $c_1/\theta\in[c_2,c_1/c_2]$ and
            \begin{equation*}
                g(\theta)=\frac{c_2-c_1}{\theta}\leq 1-\frac{c_1}{c_2}=h_{\bm{c}}\Bigl(\frac{c_1}{\theta}\Bigr).
            \end{equation*}
        \item If $\theta\in[c_1/c_2,1]$, then $c_1/\theta\in[c_1,c_2]$ and
            \begin{equation*}
                g(\theta)=(\theta-c_1)+\left(c_2-\frac{c_1}{c_2}\right)\leq \left(1-\frac{c_1}{\theta}\right)+\left(c_2-\frac{c_1}{c_2}\right)=h_{\bm{c}}\Bigl(\frac{c_1}{\theta}\Bigr).
            \end{equation*}
    \end{enumerate}
    This treats all the cases $0<\lambda<\theta<1$, as required.
\end{proof}
\subsection{Closure under suprema}\label{ss:supr}
In this section, we prove that $\mathcal{A}_d$ is closed under taking suprema.
This essentially follows since $\mathcal{A}_d$ is uniformly Lipschitz on $[0,1-\delta]$ for any $\delta>0$.
\begin{proposition}\label{p:sup-closure}
    Let $(\varphi_j)_{j\in\mathcal{J}}$ be some family of elements in $\mathcal{A}_d$.
    Then $\sup_{j\in\mathcal{J}}\varphi_j\in\mathcal{A}_d$.
\end{proposition}
\begin{proof}
    Let $f=\sup_{j\in\mathcal{J}}\varphi_j$.
    Get a sequence $J_1\subset J_2\subset\cdots\subset\mathcal{J}$ such that each $J_n$ is finite and with
    \begin{equation*}
        f_n\coloneqq\max\{\varphi_i:i\in J_n\}
    \end{equation*}
    that $f=\lim_{n\to\infty}f_n$ pointwise.
    An easy computation shows that if $\varphi_1,\varphi_2\in\mathcal{A}_d$, then $\max\{\varphi_1,\varphi_2\}\in\mathcal{A}_d$; in particular, each $f_n\in\mathcal{A}_d$.

    We first show that $f\in C([0,1])$.
    Since $(f_n)_{n=1}^\infty$ is monotonically increasing, by the Arzelà--Ascoli Theorem, it suffices to show that $(f_n)_{n=1}^\infty$ is uniformly bounded and uniformly equicontinuous.
    Uniform boundedness is immediate, so we must verify uniform equicontinuity.

    Set $b=\lim_{n\to\infty}f_n(1)$ and let $N$ be sufficiently large so that $f_n(1)>b-\varepsilon/2$ for all $n\geq N$.
    Since $f_N$ is continuous, get $\delta>0$ so that $f_N(y)>f_N(1)-\varepsilon/2$ for all $y\in[1-\delta,1]$.
    Then $|f_n(x)-f_n(y)|\leq \varepsilon$ whenever $x,y\in[1-\delta,1]$.
    Finally, since each $f_n\in\mathcal{A}_d$, the function $f_n$ is uniformly Lipschitz on $[0,1-\delta]$ as proven in \cref{p:lip-bounds}.
    It follows that $(f_n)_{n=1}^\infty$ is uniformly equicontinuous on $[0,1]$.

    Thus $f\in C([0,1])$.
    To verify that $f\in\mathcal{A}_d$, let $0<\lambda<\theta<1$ be arbitrary.
    Then for any $\varepsilon>0$, get $n$ such that $\norm{f_n-f}_\infty\leq\varepsilon$ so that
    \begin{align*}
        (1-\lambda)f(\lambda)-(1-\theta)f(\theta)&\leq (1-\lambda)(f_n(\lambda)+\varepsilon)-(1-\theta)(f_n(\theta)-\varepsilon)\\
                                                       &\leq (\theta-\lambda)f_n(\lambda/\theta)+2\varepsilon\\
                                                       &\leq (\theta-\lambda)f(\lambda/\theta)+3\varepsilon
    \end{align*}
    for any $\varepsilon>0$, so the inequality holds.
    The lower inequality follows identically.
\end{proof}
\begin{remark}
    Note that $\mathcal{A}_d$ is not compact: for example, consider the functions $\varphi_n(\theta)=\min\{c_n/(1-\theta),1\}$ with constants $c_n>0$.
    If $\lim_{n\to\infty}c_n=0$, then $\varphi_n$ converges pointwise to the function which is $0$ on $[0,1)$ and $1$ at $1$, and hence has no uniformly convergent subsequence.
    However, a simple modification of the above proof gives that for every $\delta>0$, the restriction of $\mathcal{A}_d$ to $C([0,1-\delta])$ is compact.
\end{remark}
\subsection{Characterization of upper Assouad spectra}\label{ss:ua-char}
We conclude this section with the proof of \cref{ic:ua-char}.
We recall that the family $\mathcal{M}_d$ is defined in \cref{e:Id}.
\begin{proofref}{ic:ua-char}
    To see that \cref{im:ua-as-set} implies \cref{im:ua-fun}, if $F\subset\R^d$ has $\dimAs\theta F=\varphi(\theta)$ and $\dimuAs\theta F=\overline{\varphi}(\theta)$, by \cite[Theorem~2.1]{fhh+2019},
    \begin{equation*}
        \overline{\varphi}(\theta)=\sup_{0<\theta'\leq\theta}\varphi(\theta')
    \end{equation*}
    so by \cref{c:upper-a}, $\overline{\varphi}\in\mathcal{A}_d$ so $\theta\mapsto(1-\theta)\overline{\varphi}(\theta)$ is decreasing.
    Of course, $\overline{\varphi}$ is increasing as well.

    Next, \cref{im:ua-fun} is equivalent to saying that for each $\lambda\in(0,1)$,
    \begin{equation*}
        f_{(\overline\varphi(\lambda),\lambda)}(\theta)\leq\overline{\varphi}(\theta)\cdot(1-\theta)
    \end{equation*}
    for all $\theta\in(0,1)$, with equality at $\theta=\lambda$.
    Since $\overline{\varphi}\in\mathcal{A}_d$ by \cref{l:incr-bound}, $\overline{\varphi}$ is uniformly continuous on $(0,1)$ and therefore $\overline{\varphi}=\sup_{\lambda\in\mathcal{L}}f_{(\overline\varphi(\lambda),\lambda)}$ for any countable dense subset $\mathcal{L}\subset(0,1)$.
    This implies \cref{im:ua-sup}.

    Finally, to see \cref{im:ua-sup} implies \cref{im:ua-as-set}, since $\mathcal{A}_d$ is closed under suprema by \cref{p:sup-closure}, if $\overline{\varphi}(\theta)=\sup_{f\in\mathcal{F}}f(\theta)$ for some $\mathcal{F}\subset\mathcal{M}_d$, then $\overline{\varphi}\in\mathcal{A}_d$.
    Thus the result follows by \cref{t:moran-exist}.
\end{proofref}

\section{Examples with exceptional Assouad spectra}\label{s:exceptional}

\subsection{Hölder failure at 1}\label{ss:hol-failure}
Here, we prove the following result which states that there is no control of the rate at which Assouad spectrum approaches the quasi-Assouad dimension.
\begin{proposition}
    Let $f\colon[0,1]\to[0,d]$ be a continuous increasing function with $f(0)>0$.
    Then there exists a compact set $F\subset\R^d$ such that $\dimAs\theta F\leq f(\theta)$ for all $\theta\in(0,1)$ and $\lim_{\theta\to 1}\dimAs\theta F=f(1)$.
\end{proposition}
\begin{proof}
    For $0\leq \theta<1$, let $h(\theta)=\min_{0\leq \theta'\leq\theta}(1-\theta')f(\theta')$ and let $\varphi(\theta)=h(\theta)/(1-\theta)$.
    By definition, $h$ is decreasing, $\varphi\leq f$, and since $f(\theta)\leq d$, $\lim_{\theta\to 1}h(\theta)=0$.
    Next, let us verify that $\varphi$ is increasing.
    Let $0\leq \lambda < \theta<1$ be arbitrary.
    If $h(\lambda)=h(\theta)$, then it follows immediately that $\varphi(\theta)>\varphi(\lambda)$.
    Otherwise, let $\theta'$ attain the minimum in the definition of $h(\theta)$.
    Since $h(\lambda)>h(\theta)$, we must have $\lambda\leq\theta'\leq \theta$ so that
    \begin{equation*}
        \varphi(\lambda)\leq f(\lambda) \leq f(\theta') \leq \frac{h(\theta)}{1-\theta'} \leq \varphi(\theta).
    \end{equation*}
    Therefore $\varphi\in\mathcal{A}_d$ and $\varphi(\theta)\leq f(\theta)$.

    Finally, we verify that $\lim_{\theta\to 1}\varphi(\theta)=f(1)$.
    Since $(1-\theta)f(\theta)>0$ for all $0\leq \theta < 1$ and $\lim_{\theta\to 1}h(\theta)=0$, for all $0\leq \lambda<1$, there is a $\lambda\leq\theta<1$ so that $h(\theta)=(1-\theta)f(\theta)$.
    Thus $\varphi(\theta)=f(\theta)$ for a sequence of $\theta$ converging to $1$; but $f$ is continuous so $\lim_{\theta\to 1}f(\theta)=f(1)$.

    To conclude, \cref{t:moran-exist} gives a compact set $F\subset\R^d$ such that $\dimAs\theta F=\varphi(\theta)$.
    By the properties of $\varphi$ established above, the claim follows.
\end{proof}
We also consider an explicit example.
Consider the function $f(\theta)=1+\frac{1}{\log(1-\theta)}$; note that $f$ is not Hölder at $1$.
A direct computation shows that there is some minimal $\theta_0\in(0,1)$ so that $(1-\theta)f(\theta)$ is decreasing on $[\theta_0,1]$.
Thus if we define
\begin{equation*}
    \sigma(\theta)=\begin{cases}
        \frac{(1-\theta_0)f(\theta_0)}{1-\theta} &: 0<\theta\leq\theta_0\\
        f(\theta) &: \theta_0<\theta<1
    \end{cases}
\end{equation*}
then $\sigma$ is a continuous increasing function of $\theta$ with $(1-\theta)\sigma(\theta)$ decreasing.
Thus by \cref{l:incr-bound}, $\sigma\in\mathcal{A}_1$.
A plot of $\sigma(\theta)$ and the upper bound $\min\{(1-\theta_0)f(\theta_0)/(1-\theta),1\}$ is given in \cref{f:non-lipschitz}.
\begin{figure}[t]
    \begin{center}
        \begin{tikzpicture}[>=stealth,xscale=8,yscale=5,font=\tiny]
            \draw[->] (-0.025,0) -- (1.1,0);
            \draw[->] (0,-0.025) node[below]{$0$} -- (0,1.1);
            \draw[] (1,-0.012) node[below]{$1$} -- (1,0.012);

            \node (xlabel) at (0.5,-0.1) {$\theta$};

            \matrix [draw,below right,scale=0.8,column sep=0.4em, inner sep=0.4em, align=left] at (0.05,1) {
                \draw[thick] (0,-0.15) -- (0.71,-0.15); & \node[inner sep=0pt] {$\sigma(\theta)$}; \\
                \draw[thick, dotted] (0,-0.25) -- (0.71,-0.25); & \node[inner sep=0pt] {$\min\bigl\{\frac{(1-\theta_0)f(\theta_0)}{1-\theta},1\bigr\}$}; \\
            };

            \draw[thick, smooth, variable=\x] plot file {figures/non_lipschitz_points.txt};
            \def\tzero{0.801712}
            \def\gtzero{0.381966}
            \draw[thick, dotted, smooth, variable=\x,domain=\tzero:{1-(1-\tzero)*\gtzero}] plot ({\x}, {(1-\tzero)*\gtzero/(1-\x)});
            \draw[thick, dotted, smooth, variable=\x,domain={1-(1-\tzero)*\gtzero}:1] plot ({\x}, {1});
            \draw[dotted] (\tzero,0) node[below]{$\theta_0$} -- (\tzero,\gtzero);
        \end{tikzpicture}
    \end{center}
    \caption{Plot of a spectrum which is not Hölder at $1$, along with the general upper bound.}
    \label{f:non-lipschitz}
\end{figure}

\subsection{Non-monotonicity on any open set}\label{ss:non-mono-open}
In this section, we prove that Assouad spectra which are non-monotonic on every open subset of $(0,1)$ are dense in the set of upper Assouad spectra.
Throughout this section, we fix a non-zero increasing $\varphi\in\mathcal{A}_d$, and as usual write $\beta(\theta)=(1-\theta)\varphi(\theta)$.

We recall that the functions $h_{(\kappa,c_1,c_2)}$ for $(\kappa,c_1,c_2)\in C_d$ are defined in \cref{ss:non-mono-ex}.
Fix $0<\lambda<1$ and for $0\leq y\leq\beta(\lambda)$, define
\begin{equation*}
    c(\lambda,y)=\frac{\lambda+y/\varphi(\lambda)- 1 + \sqrt{(\lambda + y/\varphi(\lambda) - 1)^2 + 4\lambda}}{2}.
\end{equation*}
The constraint on $y$ ensures that $\bm{c}(\lambda,y)\coloneqq(\varphi(\lambda),\lambda,c(\lambda,y))\in\mathcal{C}_d$.
Note that $c(\lambda,y)$ is chosen precisely so that
\begin{equation*}
    h_{\bm{c}(\lambda,y)}(\lambda)=y.
\end{equation*}
Observe that $h_{\bm{c}(\lambda,y)}\leq\beta$ by \cref{ic:ua-char}.
We also let $L_{\lambda,y}$ denote the unique affine function passing through the point $(\lambda,y)$ with slope $-\varphi(\lambda)$.
Equivalently,
\begin{equation*}
    L_{\lambda,y}(\theta) = h_{\bm{c}(\lambda, y)}(\theta)\qquad\text{for all}\qquad \lambda \leq \theta \leq c(\lambda,y).
\end{equation*}
Note that $L_{\lambda,y}$ has unique zero $\lambda+y/\varphi(\lambda)$.

Next, we define a useful family of approximations of the function $\varphi$.
For each $1\leq t<\infty$, define the functions
\begin{equation*}
    \varphi_t(\theta) = \varphi(\theta^t)\qquad\text{and}\qquad \beta_t(\theta) = (1-\theta)\varphi_t(\theta).
\end{equation*}
This family of functions uniformly approximates $\varphi$ from below while also satisfying a key ``affine partitioning'' property \cref{im:slope-extend}.
\begin{lemma}\label{l:key-param}
    Suppose $\varphi\in\mathcal{A}_d$ is strictly increasing.
    Then the following hold.
    \begin{enumerate}[nl,r]
        \item\label{im:slope-extend} Let $t\in(1,\infty)$ and $\lambda\in(0,1)$ and let $\ell \coloneqq L_{\lambda^t,\beta(\lambda^t)}$.
            Then $\ell(\lambda) = \beta_t(\lambda)$.
            Moreover, $\ell(\theta)>\beta_t(\theta)$ for $0<\theta<\lambda$ and $\ell(\theta) < \beta_t(\theta)$ for $\lambda<\theta<1$.
        \item\label{im:affine-drop} Let $\lambda\in(0,1)$.
            Then $L_{\lambda,\beta_t(\lambda)}(\theta) < \beta_t(\theta)$ for all $\lambda<\theta<1$.
        \item\label{im:monotone} For all $1\leq t_1<t_2$ and $0<\theta<1$, we have $\varphi_{t_1}(\theta)>\varphi_{t_2}(\theta)$.
        \item\label{im:beta-decr} For all $1<t<\infty$, $\beta_t$ is strictly decreasing.
        \item\label{im:Ad-member} For all $1\leq t<\infty$, $\varphi_t$ is strictly increasing and an element of $\mathcal{A}_d$.
        \item\label{im:approx} We have $\lim_{t\to 1}\norm{\varphi_t-\varphi}_\infty=0$.
    \end{enumerate}
\end{lemma}
\begin{proof}
    Let $t\in(1,\infty)$ and $\lambda\in(0,1)$, and let $\ell$ be defined as in \cref{im:slope-extend}.
    It is a direct computation that $\ell(\lambda) = \beta_t(\lambda)$.
    Moreover, since $\varphi$ is strictly increasing, the family of lines $L_{\lambda,\beta_t(\lambda)}$ is strictly increasing in the following sense: $L_{\lambda_1,\beta_t(\lambda_1)}(\theta)<L_{\lambda_2,\beta_t(\lambda_2)}(\theta)$ for all $\lambda_1<\lambda_2$ and $0<\theta<1$.
    Thus by monotonicity of $\theta\mapsto\theta^t$ and since $L_{\lambda,\beta_t(\lambda)}(\lambda) = \beta_t(\lambda)$, this implies \cref{im:slope-extend} and \cref{im:monotone}.
    Now, \cref{im:affine-drop} follows from \cref{im:slope-extend} since $\varphi$ is strictly increasing, so $L_{\lambda,\beta_t(\lambda)}(\theta)<\ell(\theta)$ for all $\lambda<\theta<1$.

    Next, since $\varphi$ is strictly increasing, it is clear that $\varphi_t$ is strictly increasing.
    Moreover, by \cref{im:slope-extend}, $L_{\lambda,\beta_t(\lambda)}(\theta)\leq \beta_t(\theta)$ for all $\theta\geq \lambda$, and since $(1-\theta)/(1-\theta^t)$ is decreasing (resp.\ strictly decreasing for $t>1$), $\beta_t$ is also decreasing (resp.\ strictly decreasing for $t>1$).
    This yields \cref{im:beta-decr} and therefore \cref{im:Ad-member} by \cref{ic:ua-char}.

    And finally, \cref{im:approx} holds since $\theta\mapsto \theta^t$ uniformly converges to the identity map on $[0,1]$ as $t\to 1$.
\end{proof}
\begin{remark}
    There is nothing particularly special about the function $\theta\mapsto \theta^t$ for $1\leq t<\infty$.
    Take any increasing homeomorphism $\phi$ of $[0,1]$ such that $\phi(\theta)\leq\theta$ and $\theta\mapsto (1-\theta)/(1-\phi(\theta))$ is decreasing.
    Then $\varphi\circ\phi$ is an increasing element of $\mathcal{A}_d$ which satisfies \cref{im:slope-extend}.
\end{remark}
We now introduce the key property for our inductive construction.
\begin{definition}
    Let $\mathcal{L}\subset(0,1)$ be a finite set, let $y\colon\mathcal{L}\to\R$, and let $1<t<2$.
    We say that the triple $(\mathcal{L},y,t)$ is \emph{zigzagging} if the following conditions hold:
    \begin{enumerate}[nl,a]
        \item The function $y$ is strictly decreasing.
        \item For all $\lambda\in\mathcal{L}$, $\beta_t(\lambda)<y(\lambda)<\beta(\lambda)$.
        \item For all $\lambda\in\mathcal{L}$ and $\theta\in\mathcal{L}\setminus\{\lambda\}$, $h_{\bm{c}(\theta,y(\theta))}(\lambda)<y(\lambda)$.
    \end{enumerate}
\end{definition}
A depiction of a zigzagging family can be found in \cref{f:zigzag}.
\begin{figure}[t]
    \begin{center}
        \begin{tikzpicture}[xscale=11, yscale=5,>=stealth]
            \draw[thick,->] (-0.05,0) -- (1.05,0);
            \draw[thick,->] (0,-0.05) node[below]{$0$} -- (0,1.05);
            \draw[thick] (1,-0.05) node[below]{$1$} -- (1,0.05);

            \begin{scope}[thick, smooth, variable=\x, domain=0:1]
                \draw[dashed] plot ({\x}, {1-(\x)^2});
                \draw[dotted] plot ({\x}, {(1 + (\x)^2)*(1 - \x)});
            \end{scope}

            \begin{scope}[thick]
                \draw (0., 0.867188) node[left]{$y(1/4)$} -- (0.25, 0.867188) -- (0.472665, 0.588856) -- (0.528915, 0.588856) -- (1., 0.);
                \draw (0., 0.6875) node[left]{$y(1/2)$} -- (0.5, 0.6875) -- (0.68658, 0.40763) -- (0.728247, 0.40763) -- (1., 0.);
            \end{scope}

            \begin{scope}[gray, dotted]
                \draw (0.25,-0.05) node[below, black]{$\frac{1}{4}$} -- (0.25, 0.867188);
                \draw (0.5,-0.05) node[below, black]{$\frac{1}{2}$} -- (0.5, 0.6875);
            \end{scope}

            \matrix [draw,below left,scale=0.8,column sep=0.3em, inner sep=0.3em, align=left] at (current bounding box.north east) {
                \draw[thick, dashed] (0,-0.2) -- (0.71,-0.2); & \node[inner sep=2pt] {$\beta$}; \\
                \draw[thick, dotted] (0,-0.2) -- (0.71,-0.2); & \node[inner sep=2pt] {$\beta_t$}; \\
                \draw[thick] (0,-0.2) -- (0.71,-0.2); & \node[inner sep=2pt] {$h_{\bm{c}(\lambda, y(\lambda)}$}; \\
            };
        \end{tikzpicture}
    \end{center}
    \caption{A plot of a zigzagging family for $\varphi(\theta)=1+\theta$, with $\mathcal{L}=\{1/4,1/2\}$, $t=2$, and an appropriate choice of $y$.}
    \label{f:zigzag}
\end{figure}
If $(\mathcal{L},y,t)$ is zigzagging, we define the corresponding function
\begin{equation}\label{e:psi-def}
    \psi = \psi_{\mathcal{L},y,t} = \max\Bigl\{\max_{\lambda\in\mathcal{L}}h_{\bm{c}(\lambda, y(\lambda))}, \beta_t\Bigr\}
\end{equation}
Observe that $\psi\leq\beta$, and moreover $\psi=h_{\bm{c}(\lambda, y(\lambda))}$ in a neighbourhood of each $\lambda\in\mathcal{L}$.

The key observation is that zigzagging families can be extended by arbitrary elements not in $\mathcal{L}$ in a way which only changes the definition of $\psi$ locally.
\begin{lemma}\label{l:mono-ext}
    Let $\varphi\in\mathcal{A}_d$ be strictly increasing and let $(\mathcal{L},y,t)$ be zigzagging with corresponding function $\psi$.
    Let $\zeta\in(0,1)\setminus\mathcal{L}$.
    Then for all $\delta>0$, there exists an extension $y(\zeta)\in\R$ such that $(\mathcal{L}\cup\{\zeta\},y,t)$ is zigzagging, $y(\zeta)\leq\psi(\zeta)+\delta$, and
    \begin{equation*}
        h_{\bm{c}(\zeta, y(\zeta))}(\theta) \leq \psi(\theta)\quad\text{for all}\quad \theta\geq \lambda+\delta.
    \end{equation*}
\end{lemma}
\begin{proof}
    Let $(\mathcal{L},y,t)$ be zigzagging with corresponding function $\psi$.
    Let $\delta>0$ be fixed.
    The proof will follow from two key observations.
    \begin{enumerate}
        \item\label{im:2-bound} \emph{For all $0\leq b\leq\beta(\zeta)$ and $\theta \geq c(\zeta,b)$, $h_{\bm{c}(\zeta,b)}(\theta) < \beta_t(\theta)$.}
            Recall that $h_{\bm{c}(\zeta,b)}(\theta)=f_{(\zeta, \varphi(\zeta))}(\theta)$ for all $\theta\geq \zeta/c(\zeta,b)\geq \zeta^{1/2}$.
            By \cref{im:slope-extend} of \cref{l:key-param}, $\beta_2(\theta)\geq f_{(\zeta,\varphi(\zeta))}(\theta)$ for all $\theta\geq\zeta^{1/2}$.
            Moreover, $\beta_t>\beta_2$ by \cref{l:key-param}~\cref{im:monotone}.
            Then the claim follows since $\beta_t$ is decreasing and $h_{\bm{c}(\zeta,b)}$ is constant on the interval $[c(\zeta,b),\zeta/c(\zeta,b)]$.

        \item\label{im:strict} \emph{We have $L_{\zeta,\psi(\zeta)}(\theta)< \psi(\theta)$ for all $\zeta<\theta<1$.}
            There are three cases.

            First, if $\psi(\zeta)=\beta_t(\zeta)$, then $L_{\zeta,\psi(\zeta)}(\theta) < \beta_t(\theta)\leq \psi(\theta)$ for all $\zeta<\theta<1$ by \cref{l:key-param}~\cref{im:affine-drop}.

            Second, suppose $\psi(\zeta) = h_{\bm{c}(\lambda, y(\lambda))}$ for some $\lambda>\zeta$.
            Since $\varphi$ is strictly increasing and $y(\lambda)<\beta(\lambda)=(1-\lambda)\varphi(\lambda)$,
            \begin{equation*}
                y(\lambda)\left(\frac{1}{\varphi(\zeta)}-\frac{1}{\varphi(\lambda)}\right) < (1-\lambda)\frac{\varphi(\lambda)}{\varphi(\zeta)} - (1-\lambda)\leq \lambda-\zeta.
            \end{equation*}
            The second inequality follows since $\beta$ is decreasing.
            Rearranging,
            \begin{equation*}
                \zeta + \frac{y(\lambda)}{\varphi(\zeta)} < \lambda + \frac{y(\lambda)}{\varphi(\lambda)}.
            \end{equation*}
            But the left hand side is the unique zero of $L_{\zeta, y(\lambda)}$ and the right hand side is the unique zero of $L_{\lambda, y(\lambda)}$.
            Thus for all $\theta>\zeta$, recalling that $\psi(\zeta)=y(\lambda)$ by assumption,
            \begin{equation*}
                L_{\zeta, \psi(\zeta)}(\theta) < h_{\bm{c}(\lambda, y(\lambda))}(\theta)\leq \psi(\theta).
            \end{equation*}

            Finally, suppose $\psi(\zeta) = h_{\bm{c}(\lambda, y(\lambda))}$ for some $\lambda>\zeta$.
            But $\varphi(\lambda)<\varphi(\zeta)$, so for all $\theta>\zeta$, since $h_{\bm{c}(\lambda, y(\lambda))}$ has slope either $0$ or $-\varphi(\lambda)$,
            \begin{equation*}
                L_{\zeta,\psi(\zeta)}(\theta)<h_{\bm{c}(\lambda, y(\lambda))}\leq \psi(\theta).
            \end{equation*}
            This treats all possible cases, as required.
    \end{enumerate}
    Now let $\delta>0$ be arbitrary.
    We may assume $\zeta+\delta<\lambda$ for all $\lambda\in\mathcal{L}$ with $\lambda>\zeta$.
    Let $s=\min\{y(\lambda):\lambda<\zeta\}$; note that $s>\psi(\zeta)$ since $(\mathcal{L},y,t)$ is zigzagging.
    By \cref{im:strict} and since $\psi$ is continuous, by choosing $\psi(\zeta)< y(\zeta)\leq\min\{s, \psi(\zeta)+\delta\}$ sufficiently small, we may assume that $L_{\zeta,y(\zeta)}(\theta) \leq\psi(\theta)$ for all $\zeta+\delta\leq \theta\leq c(\zeta, y(\zeta))$.
    But
    \begin{equation*}
        h_{\bm{c}(\zeta, y(\zeta))}(\theta)\leq\beta_2(\theta)<\beta_t(\theta)\leq\psi(\theta)
    \end{equation*}
    for all $\theta\geq c(\zeta, y(\zeta))$ by \cref{im:2-bound}.
    Thus choosing $y(\zeta)>\psi(\zeta)$ sufficiently small, the claim follows.
\end{proof}
With this result, the proof of \cref{it:non-mono} follows without too much more trouble.
Actually, we prove the following slightly more general result.
\begin{theorem}\label{t:non-mono}
    Let $\varphi\in\mathcal{A}_d$ be increasing and let $\mathcal{L}\subset(0,1)$ be an arbitrary countable set.
    Then for any $\varepsilon>0$, there exists $F\subset\R^d$ such that $f(\theta)=\dimAs\theta F$ satisfies $\norm{f-\varphi}_\infty\leq\varepsilon$ and for all $\lambda\in\mathcal{L}$, $\diniu-f(\lambda)>0$ and $\diniu+f(\lambda)<0$.
\end{theorem}
\begin{proof}
    First, since every increasing function of $\mathcal{A}_d$ can be uniformly approximated by a strictly increasing function in $\mathcal{A}_d$, we may assume that $\varphi$ is strictly increasing.
    Next, by \cref{l:key-param}~\cref{im:approx}, we may choose $1<t< 2$ sufficiently small so that $\norm{\varphi-\varphi_t}_\infty<\varepsilon$.

    Now, enumerate $\mathcal{L}=\{\lambda_n:n\in\N\}$, and write $\mathcal{L}_n=\{\lambda_1,\ldots,\lambda_n\}$.
    We inductively define a function $y\colon\mathcal{L}\to\R$ and a decreasing sequence of functions $\gamma_n$ such that for each $n\in\N$, setting $\psi_n = \psi_{\mathcal{L}_n,y,t}$ as in \cref{e:psi-def}, the following hold:
    \begin{enumerate}[nl,r]
        \item\label{im:mt} The triple $(\mathcal{L}_n, y, t)$ is zigzagging.
        \item\label{im:prop} The functions $\gamma_n$ are continuous and decreasing.
        \item\label{im:bd} $\gamma_n(\theta)>\psi_n(\theta)$ for all $\theta\in(0,1)\setminus\mathcal{L}_n$, and for $\lambda\in\mathcal{L}_n$, $\gamma_n(\lambda)=\psi_n(\lambda)$ and
            \begin{equation*}
                \diniu\pm\gamma_n(\lambda) = \diniu\pm h_{\bm{c}(\lambda, y(\lambda))}(\lambda).
            \end{equation*}
    \end{enumerate}

    We begin by choosing $y(\lambda_1)\in(\beta_t(\lambda_1),\beta(\lambda_1))$ arbitrarily; it is clear that \cref{im:mt} holds and that a function $\gamma_1$ satisfying \cref{im:prop} and \cref{im:bd} exists.

    Now suppose we have defined $y$ on $\mathcal{L}_n$ and a function $\gamma_n$ such that \cref{im:prop} and \cref{im:bd} hold.
    Let $\delta>0$ be sufficiently small such that $E_\delta\coloneqq [\lambda_{n+1}-\delta,\lambda_{n+1}+\delta]\cap\mathcal{L}_n=\varnothing$, $E_\delta\subset(0,1)$, and $\psi_n<\min\{y(\lambda):\lambda<\lambda_{n+1};\,\lambda\in\mathcal{L}_n\}$ on $E_\delta$.
    This choice is possible since $\diniu+\psi_n(\lambda)<0$ for $\lambda\in\mathcal{L}_n$ since $(\mathcal{L}_n,y,t)$ is zigzagging.
    Reducing $\delta$ more if necessary, we may also assume that $\gamma_n> \psi_n(\lambda_{n+1})+\delta$ on $E_\delta$.
    Applying \cref{l:mono-ext} with this choice of $\delta$, get a value $y(\lambda_{n+1})$ such that the triple $(\mathcal{L}_{n+1},y,t)$ is zigzagging.
    Moreover, since $\gamma_n$ is decreasing, it follows that \cref{im:bd} holds with $\psi_{n+1}$ in place of $\psi_n$.
    Therefore we may choose $\gamma_{n+1}\leq \gamma_n$ such that \cref{im:prop} and \cref{im:bd} hold.

    Finally, let $f_n(\theta)=\psi_n(\theta)/(1-\theta)$ and let $f = \lim_{n\to\infty}f_n=\sup_{n\in\N}f_n$.
    Then $f\in\mathcal{A}_d$ by \cref{l:key-param}~\cref{im:Ad-member} and \cref{p:sup-closure}, and moreover $\norm{f-\varphi}_\infty<\varepsilon$ by choice of $t$ since $\beta_t\leq \psi\leq\beta$ by construction.
    Moreover for each $\lambda\in\mathcal{D}$, by properties of the inductive construction,
    \begin{enumerate}[nl]
        \item $\psi(\lambda)=h_{\bm{c}(\lambda,y(\lambda))}(\lambda)=y(\lambda)<\beta(\lambda)$,
        \item $\diniu- \psi(\lambda)=\diniu- h_{\bm{c}(\lambda,y(\lambda))}$, and
        \item $\diniu+ \psi(\lambda)=\diniu+ h_{\bm{c}(\lambda,y(\lambda))}$.
    \end{enumerate}
    In particular, for all $\lambda\in\mathcal{D}$, since $y(\lambda)<\beta(\lambda)$,
    \begin{equation*}
        \diniu-\Bigl(\frac{h_{\bm{c}(\lambda,y(\lambda))}(\lambda)}{1-\lambda}\Bigr)>0\qquad\text{and}\qquad\diniu+\Bigl(\frac{h_{\bm{c}(\lambda,y(\lambda))}(\lambda)}{1-\lambda}\Bigr)<0.
    \end{equation*}
    Finally, the existence of the corresponding set follows by \cref{t:moran-exist}.
\end{proof}
The construction is quite flexible: the countable set $\mathcal{L}$ can be chosen arbitrarily and moreover the function $y$ at each step of the inductive construction can be chosen from an open set of parameters.
This motivates the following question.
\begin{question}
    Are ``typical'' elements of $\mathcal{A}_d$ non-monotonic?
    Does the set of functions $\varphi\in\mathcal{A}_d$ where $\varphi$ is non-monotonic on every open subset of $(0,1)$ form a residual subset of $\mathcal{A}_d$?
\end{question}
\begin{acknowledgements}
    The author thanks Amlan Banaji, Kenneth Falconer, and Jonathan Fraser for valuable comments on a draft of this paper.
    He also thanks Haipeng Chen for pointing out a mistake in an earlier version of \cref{p:lip-bounds}, and Balázs Bárány and Lars Olsen for many helpful comments which have helped to improve the exposition in this paper.
    The author was supported by EPSRC Grant EP/V520123/1 as well as the National Sciences and Engineering Research Council of Canada.
\end{acknowledgements}

@preprint{arxiv:2208.02774,
  author = {Roos, Joris and Seeger, Andreas and Srivastava, Rajula},
  eprint = {2208.02774},
  eprinttype = {arxiv},
  journal = {Studia Math.},
  title = {Spherical maximal operators on Heisenberg groups: Restricted dilation sets},
  year = {2022}
}

@thesis{ass1977,
  address = {Orsay},
  author = {Assouad, Patrice},
  eprint = {0396.46035},
  eprinttype = {zbl},
  institution = {Univ. Paris XI},
  language = {French},
  series = {Publ. Math. Orsay},
  title = {Espaces métriques, plongements, facteurs},
  type = {Thèse de doctorat d’État},
  year = {1977}
}

@article{bf2023,
  author = {Banaji, Amlan and Fraser, Jonathan M.},
  eprint = {07662347},
  eprinttype = {zbl},
  journal = {Trans. Amer. Math. Soc.},
  language = {English},
  number = {4},
  pages = {2449--2479},
  publisher = {American Mathematical Society (AMS)},
  title = {Intermediate dimensions of infinitely generated attractors},
  volume = {376},
  year = {2023}
}

@article{bff2022,
  author = {Burrell, Stuart A. and Falconer, Kenneth J. and Fraser, Jonathan M.},
  eprint = {1510.28006},
  eprinttype = {zbl},
  journal = {Monatsh. Math.},
  language = {English},
  number = {1},
  pages = {1--22},
  publisher = {Springer Science and Business Media LLC},
  title = {The fractal structure of elliptical polynomial spirals},
  volume = {199},
  year = {2022}
}

@book{bru1994,
  author = {Bruckner, Andrew M.},
  eprint = {0796.26004},
  eprinttype = {zbl},
  language = {English},
  publisher = {Providence, RI: American Mathematical Society},
  title = {Differentiation of real functions},
  volume = {5},
  year = {1994}
}

@article{ct2023,
  author = {Chrontsios Garitsis, Efstathios K. and Tyson, Jeremy T.},
  eprint = {07731255},
  eprinttype = {zbl},
  journal = {Bull. Lond. Math. Soc.},
  language = {English},
  number = {1},
  pages = {282--307},
  publisher = {Wiley},
  title = {Quasiconformal distortion of the Assouad spectrum and classification of polynomial spirals},
  volume = {55},
  year = {2022}
}

@article{fhh+2019,
  author = {Fraser, Jonathan M. and Hare, Kathryn E. and Hare, Kevin G. and Troscheit, Sascha and Yu, Han},
  eprint = {1410.28008},
  eprinttype = {zbl},
  journal = {Ann. Acad. Sci. Fenn. Math.},
  language = {English},
  number = {1},
  pages = {379--387},
  publisher = {Finnish Academy of Science and Letters},
  title = {The Assouad spectrum and the quasi-Assouad dimension: a tale of two spectra},
  volume = {44},
  year = {2019}
}

@book{fra2021,
  address = {Cambridge},
  author = {Fraser, Jonathan M.},
  eprint = {1467.28001},
  eprinttype = {zbl},
  language = {English},
  publisher = {Cambridge University Press},
  title = {Assouad dimension and fractal geometry},
  volume = {222},
  year = {2020}
}

@article{fs2023,
  author = {Fraser, Jonathan M. and Stuart, Liam},
  eprint = {1520.28005},
  eprinttype = {zbl},
  journal = {Geom. Dedicata},
  language = {English},
  number = {1},
  pages = {32},
  publisher = {Springer Science and Business Media LLC},
  title = {The Assouad spectrum of Kleinian limit sets and Patterson--Sullivan measure},
  volume = {217},
  year = {2023}
}

@article{ft2021,
  author = {Fraser, Jonathan M. and Troscheit, Sascha},
  eprint = {1479.28005},
  eprinttype = {zbl},
  journal = {Ergodic Theory Dynam. Systems},
  language = {English},
  number = {10},
  pages = {2927--2945},
  publisher = {Cambridge University Press (CUP)},
  title = {The Assouad spectrum of random self-affine carpets},
  volume = {41},
  year = {2021}
}

@article{fy2018a,
  author = {Fraser, Jonathan M. and Yu, Han},
  eprint = {1407.28002},
  eprinttype = {zbl},
  journal = {Indiana Univ. Math. J.},
  language = {English},
  number = {5},
  pages = {2005--2043},
  publisher = {Indiana University Mathematics Journal},
  title = {Assouad-type spectra for some fractal families},
  volume = {67},
  year = {2018}
}

@article{fy2018b,
  author = {Fraser, Jonathan M. and Yu, Han},
  eprint = {1390.28019},
  eprinttype = {zbl},
  journal = {Adv. Math.},
  language = {English},
  pages = {273--328},
  publisher = {Elsevier BV},
  title = {New dimension spectra: finer information on scaling and homogeneity},
  volume = {329},
  year = {2018}
}

@article{lx2016,
  author = {Lü, Fan and Xi, Li-Feng},
  eprint = {1345.28019},
  eprinttype = {zbl},
  journal = {J. Fractal Geom.},
  language = {English},
  number = {2},
  pages = {187--215},
  publisher = {European Mathematical Society - EMS - Publishing House GmbH},
  title = {Quasi-Assouad dimension of fractals},
  volume = {3},
  year = {2016}
}

@book{mt2010,
  address = {Providence, RI},
  author = {Mackay, John M. and Tyson, Jeremy T.},
  eprint = {1201.30002},
  eprinttype = {zbl},
  language = {English},
  publisher = {American Mathematical Society},
  title = {Conformal dimension},
  volume = {54},
  year = {2010}
}

@book{rob2010,
  author = {Robinson, James C.},
  eprint = {1222.37004},
  eprinttype = {zbl},
  language = {English},
  publisher = {Cambridge: Cambridge University Press},
  title = {Dimensions, embeddings, and attractors},
  volume = {186},
  year = {2011}
}

@article{zbl:07732556,
  author = {Roos, Joris and Seeger, Andreas},
  eprint = {07732556},
  eprinttype = {zbl},
  journal = {Amer. J. Math.},
  language = {English},
  number = {4},
  pages = {1077--1110},
  publisher = {Project MUSE},
  title = {Spherical maximal functions and fractal dimensions of dilation sets},
  volume = {145},
  year = {2023}
}

@article{zbl:1448.28009,
  author = {Falconer, Kenneth J. and Fraser, Jonathan M. and Kempton, Tom},
  eprint = {1448.28009},
  eprinttype = {zbl},
  journal = {Math. Z.},
  language = {English},
  number = {1-2},
  pages = {813--830},
  publisher = {Springer Science and Business Media LLC},
  title = {Intermediate dimensions},
  volume = {296},
  year = {2020}
}

@article{zbl:1509.28005,
  author = {Banaji, Amlan and Rutar, Alex},
  eprint = {1509.28005},
  eprinttype = {zbl},
  journal = {Ann. Fenn. Math.},
  language = {English},
  number = {2},
  pages = {939--960},
  publisher = {Finnish Mathematical Society},
  title = {Attainable forms of intermediate dimensions},
  volume = {47},
  year = {2022}
}
\end{document}